\documentclass[12pt]{amsart}
\usepackage{amsmath}
\usepackage{amsthm}
\usepackage{amsfonts}
\usepackage{amscd}
\usepackage{amssymb}
\usepackage{mathrsfs}
\usepackage{graphicx}
\usepackage{psfrag}

\def\ca{{\mathcal A}}
\def\cb{{\mathcal B}}
\def\cc{{\mathcal C}}
\def\cd{{\mathcal D}}
\def\ce{{\mathcal E}}
\def\cf{{\mathcal F}}
\def\cg{{\mathcal G}}

\def\cl{{\mathcal L}}

\def\cn{{\mathcal N}}

\def\cp{{\mathcal P}}

\def\car{{\mathcal R}}
\def\cs{{\mathcal S}}
\def\ct{{\mathcal T}}

\def\cv{{\mathcal V}}

\def\bc{{\mathbb C}}
\def\baf{{\mathbb F}}

\def\bn{{\mathbb N}}
\def\br{{\mathbb R}}
\def\bz{{\mathbb Z}}

\def\eps{\varepsilon}

\newtheorem{Thm}{Theorem}[section]

\newtheorem{Prop}[Thm]{Proposition}
\newtheorem{Lemma}[Thm]{Lemma}

\newtheorem{Dfn}[Thm]{Definition}

\newtheorem{exmp}[Thm]{Example}

\theoremstyle{remark}
\newtheorem{rem}[Thm]{Remark} 
\newtheorem{ack}{Acknowledgements} 
\addtocounter{section}{-1}

\begin{document}

\title[Dirac operators and spectral triples for some fractals] 
{Dirac operators and spectral triples \\ 
for some fractal sets built on curves}

\author{Erik Christensen, Cristina Ivan and Michel L. Lapidus}
 \address{  
 Department of Mathematics, University of Copenhagen,   DK-2100\newline
\indent Copenhagen, Denmark}
 \email{echris@math.ku.dk}

 \address{   
Department of Mathematics, University of Hannover,  30167\newline
\indent  Hannover, Germany }
 \email{antonescu@math.uni-hannover.de}

\address{Department of Mathematics, University of California, Riverside\\
\indent  CA 92521-0135, USA}
\thanks{The research of the third author was partially supported by the US National Science Foundation
under the research grants DMS-0070497 and DMS-0707524}
 \email{lapidus@math.ucr.edu}

 \date{\today}


 \begin{abstract}
 
We construct spectral triples and, in particular, Dirac operators, for the algebra of continuous functions
on certain compact metric spaces. The triples are countable sums of triples where
each summand is based on a curve in the space. Several fractals, like
a finitely summable infinite tree and the Sierpinski gasket, fit naturally within our framework.
In these cases, we show that our spectral triples do describe
the geodesic distance  and the Minkowski
dimension as well as, more generally, the complex fractal dimensions of the space.
Furthermore, in the case of the Sierpinski gasket, the associated Dixmier-type trace 
coincides with the normalized Hausdorff measure 
of dimension $\log 3/ \log 2$.
\end{abstract}

\maketitle

\section{Introduction}

\noindent Consider a smooth compact spin Riemannian manifold M and
the Dirac operator $\partial_M$ associated with a fixed Riemannian
connexion over the spinor bundle $S$.  Let  $D$ denote the
extension of $\partial_M$ to $H$, the Hilbert space of 
square-integrable sections (or spinors) of $S$. Then  Alain Connes showed
that  the
geodesic distance, the dimension of $M$ and the Riemannian volume
 measure can all
be described via the operator $D$ and its associated {\em Dixmier trace}.
This observation makes it possible to reformulate the geometry of a
manifold in terms of a representation of the algebra of continuous functions as
operators on a Hilbert space on which an unbounded  self-adjoint
operator $D$ acts. The geometric structures are thus 
described in such a way that points in the manifold are not mentioned at all. 
It is then possible to replace the commutative algebra of coordinates
by a not necessarily commutative algebra, and a 
way of expressing geometric properties for noncommutative
algebras has been opened. \\
\indent Connes, alone or in collaborations, has extensively
developed the foregoing idea. 
(See, for example, the paper \cite{Co1}, the book \cite{Co2}
and the recent survey article \cite{CoM}, along with
the relevant references therein.)
He has shown how an unbounded Fredholm module over the algebra
$\mathcal{A}$ of coordinates of a possibly noncommutative space $X$
specifies some elements of the (quantized) differential and integral 
calculus on $X$ as well as the metric structure of $X$.
An {\em unbounded Fredholm module} $(\mathcal{H},D)$ consists of a
representation of $\mathcal{A}$ as bounded operators on a Hilbert
space $\mathcal{H}$ and an unbounded self-adjoint operator $D$ on
$\mathcal{H}$ satisfying certain axioms. This operator is then 
called a {\em Dirac operator}
and, if the representation is faithful, the triple
$(\mathcal{A},\mathcal{H}, D)$  is called a
{\em spectral triple}. \\
 \indent The introduction of the concept of a
  spectral triple for a not necessarily commutative algebra 
$\mathcal{A}$ makes it possible to assign a subalgebra of $\mathcal{A}$
as playing the role of the 
algebra of the infinitely differentiable functions, even though
expressions of the type $(f(x)-f(y))/(x-y)$ do not make sense anymore.\\
\indent With this space-free description of some elements of differential
geometry at our disposal, a new way of investigating 
the geometry of fractal sets seems to be possible.
Fractal sets are nonsmooth 
in the traditional sense, so at first it seems quite
unreasonable to try to apply the methods of differential geometry in the
study of such sets. On the other hand, any compact set $\ct$ is completely described
as the spectrum of the C*-algebra C$(\ct)$ consisting of the
continuous complex functions on it. This means that if we can apply
Connes' space-free techniques to a dense subalgebra of C$(\ct)$,
then we do get some insight into the sort of geometric structures on
the compact space $\ct\,$ that are compatible with the given
topology.\\
\indent Already, in their unpublished work \cite{CoS}, Alain Connes and Dennis 
Sullivan have developed a `quantized calculus' on the limit sets of `quasi-Fuchsian'
groups, including certain Julia sets. Their results are presented in 
\cite{Co2}, Chapter 4, Sections $3.\gamma$ and $3.\epsilon$, along with
related results of Connes on Cantor-type sets (see also \cite{Co3}). The latter are motivated in part by 
the work of Michel Lapidus and Carl Pomerance 
\cite{LP} or Helmut Maier \cite{LM} on the geometry and spectra of fractal strings and their 
associated zeta functions (now much further expanded in the theory of 
`complex fractal dimensions'  developed in the books \cite{L-vF1,L-vF2}).\\ 
\indent In particular, in certain cases, the Minkowski dimension
and the Hausdorff measure can be recovered from
operator algebraic data. Work in this direction was pursued by
Daniele Guido and Tommaso Isola in several papers, \cite{GI1, GI2, GI3}.\\
\indent Earlier, in \cite{La3}, using the results and 
methods of \cite{KiL1} and \cite{Co2} 
(including the notion of `Dixmier trace'),
Lapidus has constructed an analogue of (Riemannian) volume measures
on finitely ramified or p.c.f. self-similar fractals and related them to the notion
of `spectral dimension', which gives the asymptotic growth rate of the spectra
of Laplacians on these fractals. (See also \cite{KiL2} where these volume measures
were later precisely identified, for a class of fractals 
including the Sierpinski gasket.)   In his programmatic paper
\cite{La4}, building upon \cite{La3}, he then investigated in many different ways the possibility 
of developing a kind of {\em noncommutative fractal geometry}, which would merge 
aspects of {\em analysis on fractals} (as now presented e.g. in \cite{Ki}) 
and {\em Connes' noncommutative geometry} \cite{Co2}. (See also parts of \cite{La2} and \cite{La3}.)
Central to \cite{La4} was the proposal to construct 
suitable spectral triples that would capture 
the geometric and spectral aspects of a given self-similar fractal, including 
its metric structure.
Much remains to be done in this direction, however, but the present work,
along with some of its predecessors mentioned above and below, address these issues 
from a purely geometric point of view. 
In later work, we hope to be able to extend our construction to address
some of the more spectral aspects.\\
\indent In \cite{CI2}, Erik Christensen and Cristina
Ivan constructed  a spectral triple for the approximately
finite-dimensional (AF) C*-algebras.  The continuous functions on the
Cantor set form an AF C*-algebra since
the Cantor set is totally disconnected.  Hence,  it was quite natural
to try to apply the general results of that paper for AF C*-algebras to this 
well-known example. In this manner, they showed in \cite{CI2} once again 
how suitable noncommutative geometry may be 
to the study of the  geometry of a fractal. 
Since then, the authors of the present article
have searched for possible spectral triples associated to other
known fractals. The hope is that these triples may be 
relevant to both {\em fractal geometry} and {\em analysis on
fractals}.
We have been especially interested in the Sierpinski gasket,
a well-known nowhere differentiable planar curve, 
because of its key role in the development of harmonic analysis on fractals. (See, for example,
\cite{Ba}, \cite{Ki}, \cite{St}, \cite{Te1}.)  
We present here a spectral triple for this fractal which recovers its geometry 
very well. In particular, it captures some of the visually detectable aspects of the gasket.
  
\begin{figure}[h]
\includegraphics[scale=0.2]{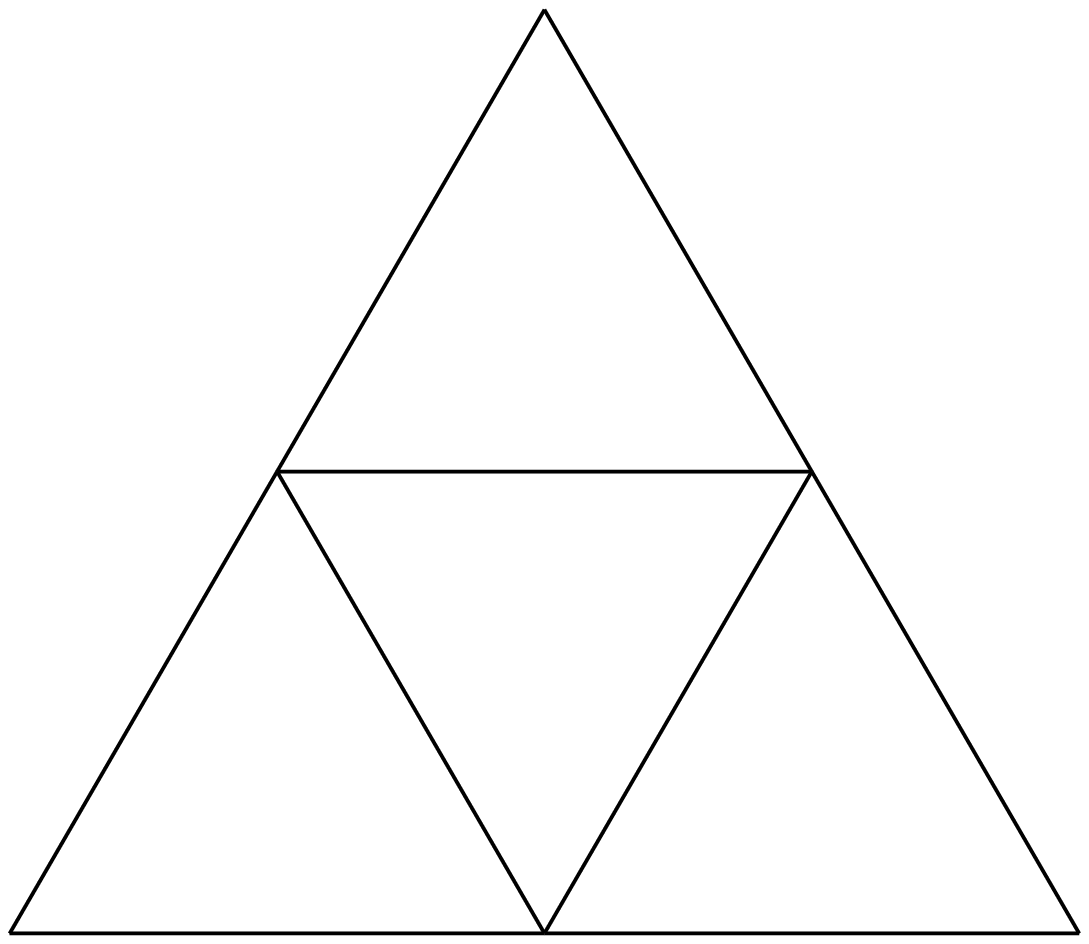}
\includegraphics[scale=0.2]{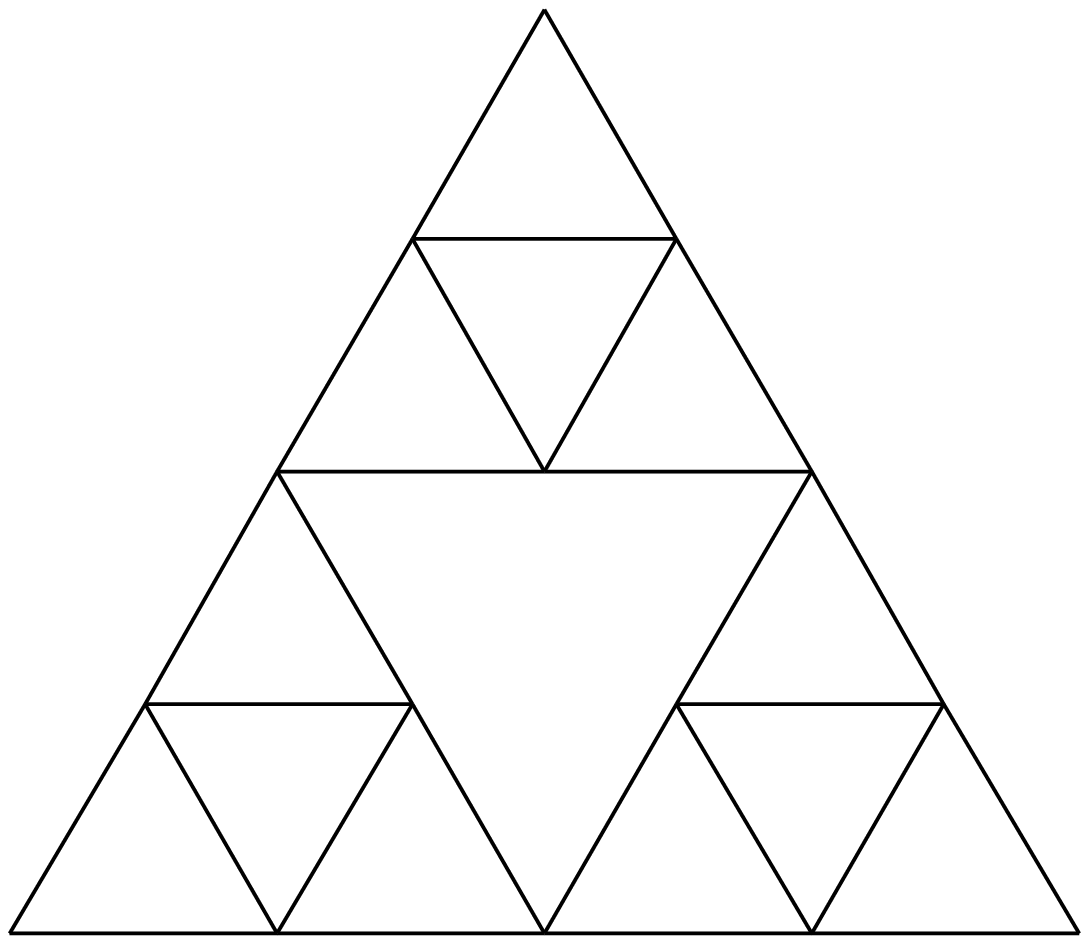}
\includegraphics[scale=0.2]{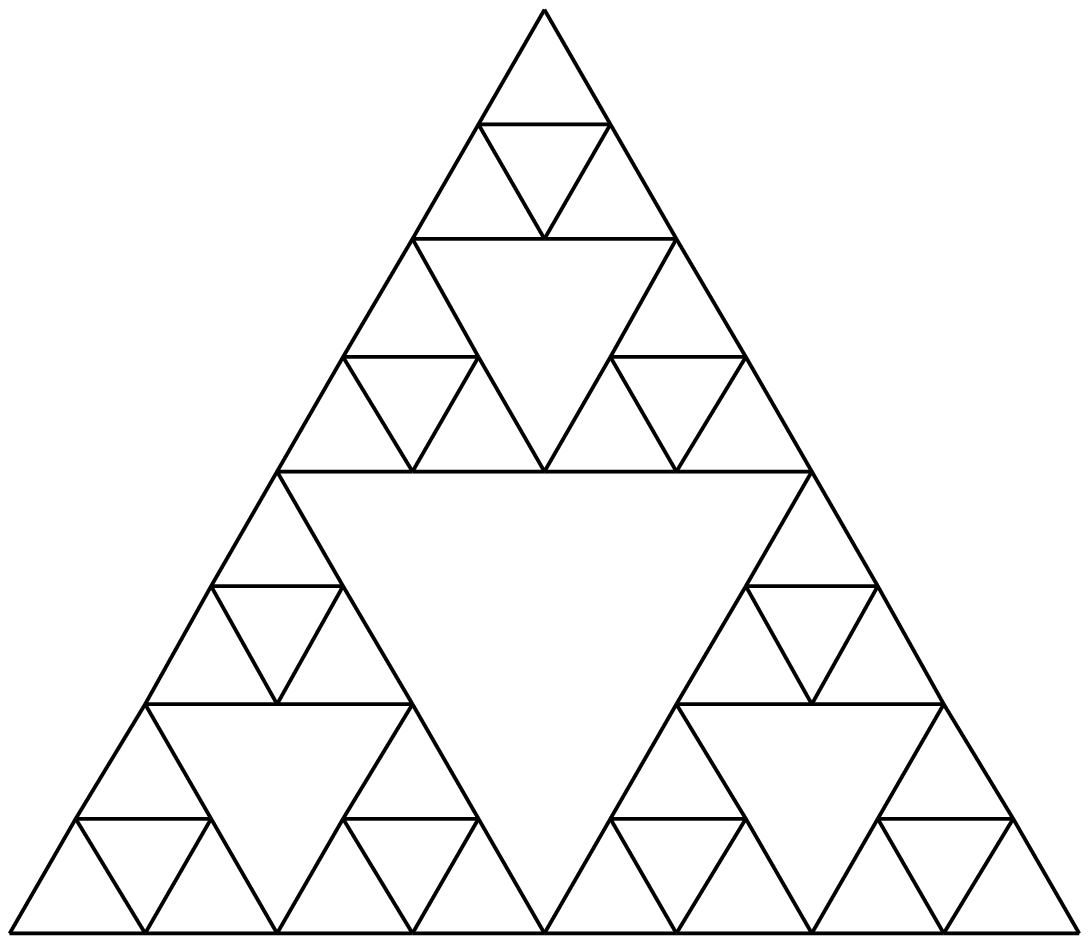}
\caption{The first 3 pre-fractal approximations to the Sierpinski gasket.}\label{SG1}
\end{figure}

\smallskip
\indent
We will now give a more detailed description of the Sierpinski gasket
and of the spectral triple we construct in that case.  
The  way of looking at the gasket, upon which our construction is
based,  is as the closure of the limit of an increasing sequence of graphs. The
pictures in Figure \ref{SG1} illustrate this.
We have chosen the starting
approximation, $\cs\cg_0$, as a triangle, and as a graph it consists of
3 vertices 
and 3 edges. The construction is of course independent of the scale,
but we will  assume that we work in the Euclidean
 plane and that the usual concept of length 
is given. We then fix the size of the largest triangle such
that all the sides have length $2\pi/3.$ 
Then,  for any natural number $n$, the
$n$-th approximation  $\cs\cg_{n} $ is constructed from
$\cs\cg_{n-1}$ by adding $3^{n}$ new vertices and $3^{n}$ new edges of
length $2^{1-n}\pi/3$. 
\smallskip
\indent
It is possible to construct a spectral triple for the gasket by
constructing a spectral triple for each new edge introduced in one of
the graphs $\cs\cg_n,$ but in order to let the spectral triples
describe the holes too, and not only the lines, 
we have chosen to look at this iterative construction as an increasing
union of triangles. Hence, $\cs\cg_0 $ consists of one equilateral triangle with
circumference (i.e., perimeter) equal to $2\pi$ and $\cs\cg_1$ consists of the union of 3
equilateral triangles each of circumference $\pi$. Continuing in this manner, we
get for any natural number $n$ that $\cs\cg_n$ is a union of $3^n$
equilateral triangles of side length $2^{1-n}\pi /3$. The advantage of this
description is that, for any equilateral triangle of circumference
$2^{1-n}\pi$, there is an obvious way  to construct a spectral triple
for the continuous functions on such a set. The idea  is to
view such a triangle as a modified circle of circumference
$2^{1-n}\pi$ and then use the standard
spectral triple for the continuous functions on such a circle. 

\indent
The direct sum of all these spectral triples is not automatically a
spectral triple, but a minor translation of each of the involved Dirac
operators can make a direct sum possible so that we can obtain a
spectral triple for the continuous functions on the Sierpinski gasket
in such a way that the following geometric structures are described
by the triple:
 
\begin{itemize}
\item[(i)] The metric induced on the Sierpinski gasket by the spectral
  triple is exactly the {\em geodesic distance} on the gasket.
\item[(ii)] The spectral triple is $s$-summable for any $s $ greater
  than \newline
  $\log 3 / \log 2$ and not $s$-summable for any $s$ smaller than or equal
  to this number. Hence it has {\em metric dimension} 
  $\log 3/\log 2$, equal to the {\em Minkowski dimension} of the gasket.  
\item[(iii)] The set of (geometric) complex fractal dimensions of the Sierpinski gasket
(defined here as the poles of the zeta function of the spectral triple)
is given by 
\begin{displaymath}
\left \{ \frac{\log 3}{\log 2}+ \sqrt{-1}\cdot k \cdot \frac{2\pi}{\log 2} \mid k\in \mathbb{Z} \right \}
\cup \{ 1 \},
\end{displaymath}
much as was suggested in \cite{La4}. 
\item[(iv)] For each natural number $n$, let the set of vertices in the
  graph $\cs\cg_n$ be
  denoted $\cv_n$ and its cardinality $\vert \cv_n \vert $. Then there exists a state $\psi$ on C$(\cs\cg)$ such
  that for any continuous complex-valued function $f$ on $\cs\cg$,
\begin{displaymath}
\underset{n \to \infty }{\lim} \frac{1}{\vert \cv_n \vert}\sum_{v \in \cv_n} f(v) 
=  \psi(f).
\end{displaymath}
This function $\psi$ is a multiple of the positive functional $\tau$ on
C$(\cs\cg)$ which the spectral triple and a Dixmier trace 
 create. Further, we show that $\psi$ coincides with the 
 {\em normalized Hausdorff measure} on the Sierpinski gasket.
\end{itemize}

\indent
The methods which we have developed for the Sierpinski gasket are based on
harmonic analysis on the circle. After having provided some necessary
background in Section 1, we will recall in Section 2 the construction of the
standard spectral triple for the continuous functions on a circle.  
Since a continuous simple image of a
circle in a compact topological space, say $X$, induces a *-homomorphism
of the continuous functions on this space onto the continuous
functions on the circle, such a closed curve in a compact space will
induce an unbounded Fredholm module over the algebra of continuous functions on
$X$.

\indent
The results for continuous images of circles in the space $X$ can be
extended to simple, i.e. not self-intersecting, 
 continuous images of intervals. We call such an image  a simple 
 continuous curve in $X$. 
The method is based on the standard embedding of
an interval in a circle, by forming a circle from an interval as the
union of the interval and a copy of the same interval and then gluing
the endpoints. This is done formally in Sections 3 and 4. We can then
construct spectral triples for C$(X)$ which are formed as direct sums
of unbounded Fredholm modules  associated to some curves in the space. 
This construction of spectral triples built on curves (and loops)
is given in Section 5 and makes good sense and is further refined for finite graphs
(Section 6), as well as for suitable infinite graphs such as trees (Section 7). In the latter section,
we illustrate our results by considering the fractal tree (Cayley graph) associated to $\mathbb{F}_2$,
the noncommutative free group of two generators. In Section 8, we then apply our results and study in detail
their consequences for the important example of the Sierpinski gasket, as described above.
Finally, in Section 9, we close this paper by discussing several open problems and proposing directions for further research.   

\begin{ack} We are grateful to Alexander Teplyaev for constructive discussions about analysis on fractals 
and for bringing the newly developed theory of `quantum graphs' to our attention, after we had completed 
a preliminary version of this paper. \\
\indent Erik Christensen expresses his gratitude towards Gilles Pisier for generous support while he was visiting  
Texas A\&M University. 
\end{ack}

\section{Notation and preliminaries}
\noindent In the book \cite{Co2}, Connes explains how a lot of geometric
properties of a compact space  $X$  can be described via the algebra of continuous
functions C$(X)$, representations of this algebra and some unbounded
operators on a Hilbert space like differential operators. The concepts
mentioned do not rely on the commutativity of the algebra C$(X)$, so it
is possible to consider a general noncommutative C*-algebra, its
representations and some unbounded operators on a Hilbert space as noncommutative
versions of structures which in the first place have been
expressed in terms of sets and geometric properties. In particular,
the closely related notions of unbounded Fredholm module and
spectral triple are  fundamental and we will recall them here.

\begin{Dfn} \label{ufm}
Let $\ca$ be a unital C*-algebra. An {\bf unbounded Fredholm module} $(H,D)$ over $\ca$ 
consists of a Hilbert space $H$ which carries a
unital representation $\pi$ of $\ca$  and an unbounded, self-adjoint
operator $D$ on $H$ such that
\begin{itemize}
\item[(i)] the set $\{a \in \ca\,|\, [D, \pi(a)]$ is densely defined and extends 
to a bounded operator on $ H\}$ is a dense subset of $\ca$,
\item[(ii)] the operator $(I+D^2)^{-1}$ is compact.
\end{itemize} 
If, in addition, {\rm tr}$((I+D^2)^{-p/2}) < \infty $ for some positive real number $p$, then 
the unbounded Fredholm module, is said to be {\bf {\em p}-summable}, or just {\bf finitely
summable}. The number $\frak{d}_{ST}$ given by 
\begin{displaymath}
\frak{d}_{ST}: =\inf \{ p>0\mid {\rm tr} ((I+D^2)^{-p/2}) < \infty \}
\end{displaymath} 
is called the {\bf metric dimension} of the unbounded Fredholm module. 
\end{Dfn}

\begin{Dfn}\label{spectrip} Let $\ca$ be a unital C*-algebra and $(H,D)$ an unbounded Fredholm module  over $\ca$.
If the underlying representation $\pi$ is faithful, then $(\ca, H, D)$ is called a 
{\bf spectral triple}.
\end{Dfn}

\indent For notational simplicity, we will denote by $(\ca, H, D)$ either 
a spectral triple or an unbounded Fredholm module, whether or not the underlying representation is faithful.

\indent
We will refer to Connes' book \cite{Co2} on several occasions later in
this article. As a standard reference on the theory of operator
algebras, we use  the books by Kadison and Ringrose
 \cite{KaR}.

\indent The emphasis of this article is actually not very much on noncommutative
 C*-algebras and noncommutative geometry, but our
 interest here is to see to what extend the spectral triples (unbounded Fredholm modules) of noncommutative geometry
 can be used to describe compact spaces of a fractal nature.  
In the papers \cite{CI1, CI2}, one can find
 investigations of this sort for the Cantor set and for general
 compact metric spaces, respectively. The metrics of compact spaces introduced 
 via noncommutative geometric methods were extensively investigated
 by Marc Rieffel. (See, for example, the papers \cite{Ri1,Ri2,Ri3,Ri4}.)

\indent Our standard references for concepts from fractal geometry are
Ed-gar's books \cite{Edg1, Edg2},
Falconer's book \cite{Fa} and Lapidus and Frank-enhuijsen's books  \cite{L-vF1, L-vF2}.

\section{The spectral triple for a circle}
\noindent Much of the material in this section is well known (see e.g. \cite{Edw}) but
is not traditionally presented in the language of spectral triples. It is useful,
however, to discuss it here because to our knowledge,
it is only available in scattered references and not in the form we need it. In particular,
we will study in some detail the domains of definition     
of the relevant unbounded operators and derivations.

\smallskip
\indent Let $C_r$ denote the circle in the complex plane with radius $r > 0$
and centered at $0.$ As usual in noncommutative geometry, we are not
studying the circle directly but rather a subalgebra of the algebra of
continuous functions on it. From this point of view, it seems easier to look at
the algebra of complex continuous $2\pi r$-periodic functions on the
real line. We will let $\ca C_r$ denote this algebra. We will let $(1/2
\pi r )\frak{m}$ denote the normalized Lebesgue measure on the interval
$[-\pi r, \pi r]$ and let $\pi_r$ be the standard representation of
$\ca C_r$ as multiplication operators on the Hilbert space $H_r$ which is
defined by $H_r := L^2([-\pi r, \pi r], (1/2 \pi r )\frak{m})$.  The space
$H_r$ has a canonical orthonormal basis, denoted $\left ( \phi_k^r \right )_{k \in \bz }$,
which consists of functions in
$\ca C_r$ given by
\begin{displaymath} 
\forall \, k \in \, \bz , \quad \phi_k^r(x) \,:=\, \exp \left (\frac{ikx}{r} \right).
\end{displaymath} 
These functions are eigenfunctions of the differential operator
$\frac{1}{i}\frac{\mathrm{d}}{\mathrm{d}x}$ and the corresponding
eigenvalues are $\{k/r\,|\, k \in \bz\,\}$. The natural choice for
the Dirac operator for this situation is the closure of the
restriction of the above operator to the linear span of the basis
$\{\phi_k^r\,|\, k \in \bz \}$. We will let $D_r$ denote this
operator on $H_r$. It is well known that $D_r$ is self-adjoint
and that dom$(D_r)$, the domain of
definition of $D_r$, is given by 
\begin{displaymath}
\forall  f \in H_r: \quad  f  \in \mathrm{dom}D_r
 \Leftrightarrow \sum_{k \in \bz} \frac{k^2}{r^2}|\langle f 
 \, |\, \phi_k^r\rangle|^2 \, < \, \infty,
\end{displaymath}
where $\langle \cdot| \cdot \rangle $ is the inner product of $H_r$.\\ 
\indent For an element $f \in$ dom$D_r$, we have $D_r f \, = \,
\sum_{k \in \bz}(k/r)\langle f \, |\, \phi_k^r\rangle   \phi_k^r.$
The self-adjoint operator $D_r$ has spectrum $\{k/r\, | \, k \in  \bz
\,\}$ and each of its eigenvalues 
has multiplicity 1. Furthermore, any continuously differentiable  $2\pi 
r$-periodic function $f$ on $\br$ satisfies 
\begin{displaymath}
 [D_r, \pi_r(f)] \, = \, \pi_r(-if^{'}), 
\end{displaymath}
 so we obtain a spectral triple associated to the circle $C_r$
in the following manner.

\begin{Dfn}
The {\bf natural spectral triple}, $ST_n(C_r)$,  for the circle algebra
$\ca C_r$ is defined by  $ST_n(C_r)\, :=\, 
( \ca C_r, H_r,  D_r)$.   
\end{Dfn}

\smallskip
\indent
One of the main ingredients in the arguments to come is the possibility
to construct {\em interesting }  spectral triples as direct sums of
 unbounded Fredholm modules, each of which only carries very little information
 on the total space. In the case of the natural triples for circles,
 the number $0$ is always an eigenvalue and hence, if the sum
 operation is done a countable number of times, the eigenvalue $0$ will
 be of infinite multiplicity for the Dirac operator which is obtained
 via the  direct sum construction. In order to avoid this problem, we will, for
 the $C_r$-case,
 replace the Dirac operator $D_r$ by a slightly modified one,
 $D_r^t$, which is the   translate  of $D_r$ given
 by 
\begin{displaymath}
D_r^t \, : = \, D_r \, + \, \frac{1}{2r}I.
\end{displaymath}  
The set of eigenvalues now becomes $\{(2k+1)/2r\, |\, k \in \bz\,
\}$, but the domain of definition is the same as 
 for $D_r$  and, {\em in
  particular,} for any function  $f \in \ca C_r$, we have 
$[D_r^t, \pi_r(f) ] \, = \,  [D_r, \pi_r(f) ] $. Hence, the translation
does not really change the effect of the spectral triple.

\begin{Dfn} \label{transpt}
The {\bf translated  spectral triple}, $ST_t(C_r)$,  for the circle algebra
$\ca C_r$ is defined by  $ST_t(C_r)\, :=\, (\ca C_r, H_r, D_r^t)$. 
\end{Dfn}

\smallskip
\indent
The next question is to determine for which functions $f $ from $ \ca C_r$
the commutator $[D_r^t, \pi_r(f)] $ is bounded and densely
defined. This is done in the following lemma which is standard, but which
we include since its specific statement cannot be easily found 
in the way we need it. On the other hand, the proof is based on
elementary analysis and is therefore omitted.

\begin{Lemma} \label{domD}
Let $f \in \ca C_r$.  Then the following conditions are
equivalent: 
\begin{itemize}
\item[(i)] $[D_r^t, \pi_r(f)] $ is densely defined and
  bounded.
\item[(ii)] $ f \in\,$  {\rm dom}$(D_r)$ and $D_rf$ is
  essentially bounded.
\item[(iii)] There exists a measurable, essentially bounded function
  $g$   on the interval $[-\pi r , \pi r ]$ such
  that 
\begin{displaymath}
\int_{-\pi r}^{\pi r} g(t)\mathrm{d}t \, = \, 0 \text{  and  } 
\forall x \in [-\pi r, \pi r ] : f(x) = f(0) +
  \int_0^x g(t)\mathrm{d}t.
\end{displaymath}
\end{itemize}
If the conditions above are satisfied, then $g(x) \, = \, (iD_rf)(x) $
a.e.
\end{Lemma}

\indent
We will end this section by mentioning some properties of this
spectral triple. We will not prove any of the claims below, 
as they are easy to verify.
First, we remark that all the statements below hold for both of the triples
$ST_t(C_r)$ and $ST_n(C_r)$, although they are only formulated for
$ST_n(C_r)$.
\begin{Thm} \label{circledist}
Let $r>0 $ and let $(\ca_rC, H_r,
D_r)$ be the $ST_n(C_r)$ circle spectral triple. Then
the following two results hold: 
\begin{itemize}
\item[(i)] The metric, say $d_r$, induced by the spectral triple
  $ST_n(C_r)$ on the circle is the geodesic distance on $C_r$. 
\item[(ii)] The spectral triple $ST_n(C_r)$ is summable for any $s > 1$,
  but not for $s=1$. Hence, it has metric dimension $1$.

\end{itemize}
\end{Thm}

\section{The interval triple}

\noindent The standard way to study an interval by means of the  theory for the
circle is to take two copies
of the interval and then glue them together at  the endpoints.
When working with algebras instead of spaces, this construction is done
via  an injective
homomorphism $\Phi_{\alpha}$  of the continuous functions on the interval $[0,
\alpha]$ into the
continuous functions on the double interval $[-\alpha, \alpha] $ defined
by 
\begin{displaymath}
\forall f \in \mathrm{C}([0,\alpha]), \, \forall t \in [-\alpha, \alpha]
 : \quad  \Phi_{\alpha}(f)(t)\, :=\, f(|t|).
\end{displaymath}

The continuous functions on the interval $[0, \alpha ] $ are mapped by
this procedure onto the even continuous $2\alpha$-periodic
continuous functions on the real line, and the theory describing the
properties of the spectral triple $(\ca C_\frac{\alpha}{\pi},
  H_\frac{\alpha}{\pi}, D^t_\frac{\alpha}{\pi})$ may be used to describe a
  spectral triple for the algebra C$([0, \alpha])$. 
 We will start by defining the spectral triple
$ST_{\alpha }$ which we associate to the interval $[0,\alpha] $ before
we actually prove that this is a spectral triple. The arguments 
showing that $ST_{\alpha} $  is a spectral triple follow from the
analogous results for the circle.

\begin{Dfn} \label{intervaltrip}
Given $\alpha > 0$, the{ \mathversion{bold}$\alpha$}-{\bf interval spectral triple} $ST_{\alpha}\,
:=\, (\ca_{\alpha}, H_{\alpha},
D_{\alpha})$ is defined by
\begin{itemize}
\item[(i)] $\ca_{\alpha} \, =\, \mathrm{C}([0, \alpha])\, $;
\item[(ii)] $ H_{\alpha} \,= \, L^2([-\alpha, \alpha], \frak{m}/2\alpha)$, where the
  measure $\frak{m}/2\alpha$ is the normalized Lebesgue measure; 
\item[(iii)] the representation $\pi_{\alpha}: \ca_{\alpha} \to
  B(H_{\alpha})$ is defined for $f$ in $\ca_{\alpha}$ as the
  multiplication operator which multiplies by the function
  $\Phi_{\alpha}(f)$; 
\item[(iv)] an orthonormal basis $\{e_k\, |\, k \in
    \bz \} $ for $H_{\alpha}$ is given by \newline $e_k(x) := exp\, (i\pi k x/\alpha)$ and
    $ D_{\alpha} $ is the 
    self-adjoint operator on $H_{\alpha}$ which has all the vectors $e_k$ as
    eigenvectors and such that $D_{\alpha}e_k = (\pi k /\alpha) e_k$ for each $k\in \bz$.
\end{itemize}
\end{Dfn}   

\indent Let us look at the expression $[D_{\alpha},\pi_{\alpha}(f)],$  for a
smooth function $f$. Then it is well known that 
$[D_{\alpha},\pi_{\alpha}(f)] = \pi_{\alpha}( D_{\alpha}f),$ and we want
to remark that for an even function $f$ we have that $D_{\alpha}f$ is
an odd function. Hence, here in the most standard commutative example, we already meet a
noncommutative  aspect of the classical theory, in the sense that for
any function $f$ in $\ca_{\alpha}$ the commutator $[D_{\alpha},
\pi_{\alpha}(f)], $ if it exists, is no  longer
in the image of $\pi_{\alpha}(\ca_{\alpha})$.\\  
\indent The following proposition demonstrates that many even continuous
functions $ f$ do have bounded commutators with $D_{\alpha}$. The result
  follows directly from Lemma  \ref{domD}.

\begin{Prop} \label{cordomD}
Let $f$ be a continuous   real and even function on the interval
$[-\alpha , \alpha] $ such that $f$ is boundedly and continuously 
differentiable
outside a set of finitely many points. Then $f$ is in the domain of
definition of $D_{\alpha}$ and $D_{\alpha}f$ is bounded outside a set
of finitely many points.
\end{Prop}

\indent
The close connection between the spectral triples for the circle and
the interval yields immediately the following corollary of Theorem \ref{circledist}.

\begin{Thm} \label{intervdist}
Given $\alpha >0$, let $(\ca_{\alpha}, H_{\alpha},
D_{\alpha})$ be the $\alpha$-interval spectral triple. 
Then, for any pair of reals $s, t$ such that
  $ 0 \leq s < t \leq \alpha$,  we have 
\begin{displaymath} 
|\, t - s \,| = \sup\{ \, |\,f(t)  - f(s)\,| \,\,|\, \|\, [D_{\alpha},
\pi_{\alpha}(f)]\,\| \leq 1 \, \}.
\end{displaymath}
Furthermore, the triple is summable for any real $s > 1$ and not summable for $s =1$.
Hence, it has metric dimension 1.
\end{Thm}

\section{The $r$-triple, $ST_r$}

\noindent Let now $\ct$ be a compact and Hausdorff space and $r: [0, \alpha] \to \ct $ a
continuous and injective mapping. The image in $\ct$, i.e. the set  of
points  $ \car = r([0,\alpha])$, is called a continuous curve and $r$ is then
 a parameterization of $\car$. As usual, a continuous curve may
have many parameterizations and one of the possible uses of the concept
of a spectral
triple is  that it can help us to distinguish between various 
parameterizations of the same continuous curve. Below we will associate
to a parameterization $r$  of a continuous curve $\car$ an unbounded
Fredholm module, $ST_r$. After having read the statement of the following proposition, 
it will be clear how this module  may
be defined.

\begin{Prop} \label{STr}
Let  $r: [0,\alpha] \to \ct$ be a continuous injective mapping 
 and  $(\ca_{\alpha}, H_{\alpha}, D_{\alpha})$ the
$\alpha$-interval spectral triple.

\indent Consider the triple  $ST_r$ defined by  $ST_r\, := \, 
(\mathrm{C}(\ct), H_{\alpha}, D_{\alpha})$, where the
representation $\pi_r: \mathrm{C}(\ct) \to B(H_{\alpha})$ is defined
via a homomorphism $\phi_r \, \mathrm{ of \, \, C}(\ct) $ onto $\ca_{\alpha}$ as follows:
\begin{itemize} 
\item[(i)] $\forall  f \in \mathrm{C}(\ct),\, \forall s \in [0,\alpha]:
 \, \, \phi_r(f)(s)\,:=\,f(r(s));$
\item[(ii)]  $\forall  f \in \mathrm{C}(\ct),\,\, \pi_r(f) \, := \,
  \pi_{\alpha}(\phi_r(f))$.
\end{itemize}
Then $ST_r$ is an unbounded Fredholm module, which is summable for any $s > 1$ and not summable for $s =1$. 
\end{Prop}

\begin{proof}
The Hilbert space, the representation and the Dirac operator are
mostly inherited from the $\alpha$-interval triple, so this makes it
quite easy to verify  the properties demanded by an unbounded Fredholm module.
The only remaining problem is to prove that the subspace
\begin{displaymath}
LC := \{f \in \mathrm{C}(\ct)\, |\, [D_{\alpha},\pi_r(f)] \text{ is
  bounded}\}
 \end{displaymath}
is dense in C$(\ct)$.
 The proof of this may be based on the Stone--Weierstrass Theorem. We
 first remark that  the Leibniz differentiation rule implies that 
$LC$ is a unital self-adjoint algebra. 
 We then only have to prove that the functions in
 $LC $ separate the points of $\ct$. An argument proving this can be
 based on  Urysohn's Lemma and Tietze's Extension Theorem.\\ 
\end{proof}
\indent
We can then associate as follows an unbounded Fredholm module to a continuous curve.
\begin{Dfn}  \label{DefSTr}
Let  $r: [0,\alpha] \to \ct$ be a continuous injective mapping 
 and  $(\ca_{\alpha}, H_{\alpha}, D_{\alpha})$ the
$\alpha$-interval spectral triple. The unbounded Fredholm module  $ST_r\, := \, 
(\mathrm{C}(\ct), H_{\alpha}, D_{\alpha})$ is then called the {\bf unbounded Fredholm module
associated to the continuous curve } $r$.
\end{Dfn}
\indent
\indent We close this section by computing the metric $d_r$ on $\ct$ induced by the
parameterization $r$ of $\car$.

\begin{Prop} \label{dr}
Let  $r: [0,\alpha] \to \ct$ be a continuous injective mapping, 
 and   $ST_r\, = \, 
(\mathrm{C}(\ct), H_{\alpha}, D_{\alpha})$ the unbounded Fredholm module 
associated to $r$. The metric $d_r$ induced on $\ct$ by $ST_r$ is
given by
\begin{displaymath}
 d_r(p,q) = \begin{cases}
0 & \text{ if } p\,=\, q, \\
\infty  & \text{ if } p\,\neq \, q \text{ and } ( p \notin \car \text{ or } q
\notin \car), \\
|r^{-1}(p) - r^{-1}(q)|  & \text{ if } p\,\neq \, q \text{ and }
  p \in \car \text{ and } q \in \car. 
\end{cases}
\end{displaymath}
\end{Prop}

\begin{proof}

If one of the points, say $p$,  is not a point on the curve and $q
\neq p$,  then by
Urysohn's Lemma  there exists a nonnegative continuous 
function $f$ on $\ct$
such that $f(p)=1$, $f(q)= 0$ and, for any point $r(t)$ from  $\car$, we
have $f(r(t)) =0. $ This means that $\pi_r(f) = 0$; so for any
$N \in \bn $, we have $\|\,[D_{\alpha}, \pi_r(Nf)]\,\| \leq 1$. Hence,
$d_r(p,q) \geq N$ and $d_r(p,q) = \infty.$ 

\indent
Suppose now that  both $p$ and $q$ belong to $\car$ and $p \neq q$. Then it follows that the homomorphism $\phi_r :C(\ct) \to C([0,\alpha])$,
as defined in Proposition \ref{STr}, has the property that $r$-triple
factors through the \newline $\alpha$-interval triple via the identities
presented in Proposition \ref{STr}, and we deduce that
\begin{displaymath}
\forall f \in  \mathrm{C}(\ct):\,  \|\,[D_{\alpha}, \pi_r(f)]\| \leq 1
\iff  \|\,[D_{\alpha}, \pi_{\alpha}(\phi_r(f)]\|
\leq 1.
\end{displaymath}
Hence, on the curve $r[0, \alpha])$, the metric  
induced by the spectral triple for
the curve is simply the metric which $r$ transports from the interval
$[0,\alpha] $ to the curve.
The proposition follows.
\end{proof}

\section{Sums of curve triples}

\noindent Let $\ct$ be a compact and Hausdorff space and suppose that  for $1 \leq i \leq k, $
we are given continuous curves $r_i : [0, \alpha_i] \to \ct$. It is then
fairly easy to define the direct sum of the associated $r_i$-unbounded
Fredholm modules, but the sum may fail to be an unbounded Fredholm module in the sense that
there may not exist a dense set of functions which have bounded commutators
with all the Dirac operators simultaneously. It turns out that if
the curves do not overlap except at finitely many points, then this
problem can not occur.

\begin{Prop}  \label{sumSTr}
Let $\ct$ be a compact and Hausdorff space and for $1 \leq i \leq h,$ let $ r_i : [0,
\alpha_i]$  be a continuous curve. If for each pair $i\neq j$ the number
of points in $r_i([0,\alpha_i])\cap r_j([0,\alpha_j])$ is 
finite, then the direct sum $\oplus^{h}_{i=1} ST_{r_i}$ is an unbounded Fredholm module
for C$(\ct)$.
\end{Prop}
 
\begin{proof}

\indent For each $i \in \{1, \dots , h\},$ the $r_i$-spectral triple,
$ST_{r_i}$ 
is defined by $ST_{r_i}\, := \, 
(\mathrm{C}(\ct), H_{\alpha_i}, D_{\alpha_i})$. Further, the direct sum  
is defined by 
\begin{displaymath}
\underset{i=1}{\overset{h}{\oplus}} ST_{r_i}\, := \, 
\big(\mathrm{C}(\ct),\underset{i=1}{\overset{h}{\oplus}} H_{\alpha_i}, \underset{i=1}{\overset{h}{\oplus}}
D_{\alpha_i}\big).
\end{displaymath}
In order to prove that 
\begin{displaymath}
\left \Vert \left [ \underset{i=1}{\overset{h}{\oplus}} D_{\alpha_i},
 \left (\underset{i=1}{\overset{h}{\oplus}} \pi_{\alpha_i}\right )(f) \right ] \right \Vert
< \infty 
\end{displaymath}
 for a dense set of continuous complex functions on
$\ct$, we will again appeal to the Stone--Weierstrass theorem, and repeat
most of the
arguments we used to prove that any $ST_r$ is an unbounded Fredholm module.
 We therefore
define a set of functions $\ca$ by
\begin{displaymath}
\ca \, := \, \{ f \in \mathrm{C}(\ct)\, | \, \forall i,\, 1 \leq i \leq
h, \, \|\,[D_{\alpha_i},\pi_{\alpha_i}(f)]\,\| < \infty \}.
\end{displaymath}
\indent
As before, we see that $\ca$ is a self-adjoint  unital
algebra. Moreover, we find, just as in the proof of Proposition \ref{STr},
that given two points $p$ and $q$ in $\ct$ for which at least one of them does
not belong to the union of the points on the curves, $\cup_{i=1}^{h}
\car_i,$ there is a function $f$ in $\ca$ such 
that $f(p) \neq f(q).$ Let us then suppose that both $p$ and $q$ are
points on the curves, say $p \in \car_j$ and $p = r_j(t_p) $, $q \in
\car_k$ and $q = r_k(t_q)$. We will
then define a continuous function $g$ on the compact set $\cup_{i=1}^{h} \car_i$
by first defining it on $\car_j$, then extending its definition to
all the curves which meet $\car_j$, and then to all the curves which
meet any of these, and so on. Finally, it is defined to be zero
on the rest of the curves, which are not connected to $\car_j.$
 The definition of $g$ on $\car_j$ is given in
the following way. First, we define a continuous function $c_j$ in
$\ca_{\alpha_j}$ by 
\begin{align*} 
& c_j(t_p) := 1;\\   
&\text{if } 0  \neq t_p \, , \text{ then } c_j(0): =0 ;\\
&\text{if } t_p \neq \alpha_j \, , \text{ then } c_j(\alpha_j): = 0; \\
&c_j \text{ is extended to a piecewise affine function on }  [0, \alpha_j]; \\
&c_j \text{ is even on } [-\alpha_j, \alpha_j].
\end{align*}
\indent
By Proposition \ref{cordomD}, such a function $c_j$
 is in the domain of $D_{\alpha_j} $ and
$D_{\alpha_j}c_j$ is a bounded function. We will then transport  $c_j $ to
a continuous function $g$ on $\car_j$  by
\begin{displaymath} 
\forall t \in [0, \alpha_j]: \,g(r_j(t)) \, : = \, c_j(t).
\end{displaymath}
\indent
The next steps then consist in extending $g$ to more curves by adding one
more curve at each step.  In order to do this, we will
give  the necessary induction argument. Let us then suppose that $ g$ is
already defined on the curves $ \car_{i_1}, \dots , \car_{i_n}$ and
assume that there is yet a curve $\car_{i_{n+1} } $ on which $g$  is
not defined but  has the property that $\car_{i_{n+1}} $  intersects at least one of the curves
$ \car_{i_1}, \dots , \car_{i_n}$. Clearly, there are  only
finitely many points $s_1, \dots, s_m$ on $\car_{i_{n+1}}$ which are also
in the union of the first $n$ curves. Suppose that these
points are numbered in such a way that 
\begin{displaymath}
0 \leq r^{-1}_{i_{n+1}}(s_1) < \dots < r^{-1}_{i_{n+1}}(s_m) \leq
\alpha_{i_{n+1}}. 
\end{displaymath}    
We can then define a piecewise affine, continuous and
even function $c_{i_{n+1}}$ on $[-\alpha_{i_{n+1}},
\alpha_{i_{n+1}}] $ by 
\begin{align*} 
&  \forall l \in \{1, \dots , m\}:\, 
c_{i_{n+1}}( r_{i_{n+1}}^{-1}(s_l)): =g(s_l);\\
&\text{if } 0  \neq s_1 \, , \text{ then } c_{i_{n+1}}(0): =0 ;\\
&\text{if } s_m \neq \alpha_{i_{n+1}} \, ,  \text{ then }
c_{i_{n+1}}(\alpha_{i_{n+1}}): = 0 ;\\ 
&c_{i_{n+1}} \text{ is extended to a piecewise affine function on }
[0, \alpha_{i_{n+1}}]; \\ 
&c_{i_{n+1}} \text{ is even on }  [-\alpha_{i_{n+1}},\alpha_{i_{n+1}}].
\end{align*}
\noindent This function is in the domain of $D_{\alpha_{i_{n+1}}}$ and
$D_{\alpha_{i_{n+1}}}c_{i_{n+1}}$ is bounded. The function
$c_{i_{n+1}}$
 can then be lifted to $\car_{i_{n+1}}$ and we may define
$g$ on this curve by 

\begin{displaymath} 
\forall t \in [0, \alpha_{i_{n+1}}]: \, g(r_{i_{n+1}} (t)) \, : = \,
c_{i_{n+1}}(t). 
\end{displaymath}

\indent
This process will stop after a finite number of steps and will
yield a continuous function $g$ on the union of the curves which are
connected to $p$. On the other curves, $g$ is defined to be $0$, and  by
Tietze's extension theorem, this function can then be extended to a
continuous real function on all of $\ct$. By construction, the function $g$
belongs to $\ca$. Further, its restriction to the
union $\cup \car_{i=1}^{h}$ has the property that $g(p) = 1$ and for all other
points s in $\cup \car_{i=1}^{h}$, we have $g(s) < 1$. In particular, $g(q) \neq
g(p)$. Hence, the theorem follows.

\end{proof}

\section{ Parameterized  graphs}
\smallskip
\noindent
We will now turn to the study of {\em finite graphs} where each edge, say
$e$,  is
equipped with a weight or rather a length, say $\ell(e)$. 
Such a graph,  is called a
{\em weighted graph}. 
This concept has occurred in many places and it is very well described
in the literature on graph theory and its applications. Our point of
view on this concept is related to the study of the so-called   {\em quantum
  graphs}, although it was developed independently of it. 
We refer the reader
to the survey article  \cite{Ku1} by Peter Kuchment on this kind of
weighted graphs and their properties. 
In Section 2 of the article by  Kuchment, a quantum graph is described as a
graph for which each edge $e$ is considered as a line segment connecting
two vertices.  Any edge $e$ is given a length $\ell(e)$ and then
equipped with a parameterization, as described in the definition just below; 
it then becomes homeomorphic to the
interval $[0,\ell(e)]$. This allows one to think of the quantum graph
as a one-dimensional simplicial complex. 
\begin{Dfn} \label{pargraph}
A weighted graph  $\cg$ with vertices $\cv$, 
edges $\ce$ and a length function $\ell,$ 
is said to have a {\bf parameterization} if it is 
 a subset of a compact and Hausdorff space $\ct$ such that the
vertices in $\cv$ are points in $\ct$ and for each edge $e$ there exists a pair
of vertices $(p,q)$ and 
a continuous curve 
$\{e(t) \in \ct \, |\, 0 \leq t \leq  \ell(e)\}$ without
self-intersections such that $e(0) = p$ and
$e(\ell(e))= q.  $ Two edges $e(t)$ and $f(s)$ may only intersect at
some of their
endpoints. Given $e\in \ce $, the length of the defining interval, $\ell(e)$, is
called the length of the edge $e$. The set of all points on the edges
is denoted $\cp$.  
\end{Dfn}

\indent Even though the foregoing definition indirectly points at an orientation on
each edge, this is not the intention. Below we will allow {\em traffic in both
  directions}
on any edge, but we will only have one parameterization. 
One may wonder if any finite
graph has a parameterization as a subset of some space $\br^n$, but it
is quite obvious that if one places as many points as there are
vertices inside the space $\br^3$ such that the diameter of this set
is at most half of the length of the shortest edge, then the points
can be connected as in the graph with smooth nonintersecting strings
having the correct lengths. Consequently, it makes sense to speak of the
set of points $\cp$ of the graph and of the distances between the points
in $\cp$.
\smallskip

\indent The graphs we are studying
have natural representations as compact subsets of $\br^n$, so we
would like to have this structure present when we consider the graphs.
Hence, we would like to think of our graphs as embedded in some
compact metric space. In the case of a finite graph, which we are considering  
in this section, the ambient space may be nothing but the
simplicial complex, described above, and  equipped with a
parameterization of the 1-simplices.\\
\indent  Each edge $e$ in a parameterized graph has an internal
metric given by the parameterization and then, as in \cite{Ku1}, it is possible 
to define the distance between any two points from the same 
connected component of   $\cp$ as the length of the 
shortest path connecting the points. 
In order to keep this paper self-contained, however, we will next introduce the notions of  a path and length of a path
in this context.  
We first recall the graph-theoretic
concept of a path as a set of edges $\{e_1, \dots , e_k\}$
such that $e_i \cap e_{i+1} $ is a vertex and no vertex appears more
than once in such an intersection. It is allowed that the starting
point of $e_1$ equals the endpoint of $e_k$, in which case we say
that the path is closed, or that it is a {\em cycle}. 
  We can then extend the  concept of a path 
to a parameterized graph as follows.

\begin{Dfn} \label{path}
Let $\cg = ( \cv , \ce) $ be a graph with parameterization. Let $e, f$
be two edges, $e(s)$ a point on $e$ and $ f(t) $ a point on
$f$. Let $\cg \cup\{e(s),f(t)\}$ be the parameterized graph obtained
  from $\cg$ by adding the  vertices  $e(s), f(t)$ and, accordingly,
  dividing
some of the edges into 2 or 3 edges. A {\bf path} from $e(s) $ to
$f(t)$ is then an ordinary graph-theoretic path in $\cg
\cup\{e(s),f(t)\}$ starting at $e(s) $ and ending at $f(t)$.
\end{Dfn}

\smallskip
\indent The intuitive picture is that a path from $e(s) $ to $f(t) $ starts at
$e(s)$ and runs along the edge upon which $e(s)$ lies to  an endpoint of
that edge. From there, it continues on an ordinary graph-theoretic 
 path to an endpoint
of the edge upon which $f(t)$ is placed, and then continues on this
edge to $f(t)$. This is also almost what is described in the definition
above, but not exactly. The problem which may occur is that $e(s)$ and
$f(t)$ actually are points on the same edge, and then we must allow
the possibility that the path runs directly from
$e(s)$ to $f(t)$ on this edge, as well as 
the possibility that the path starts  from $e(s)$, runs away from $f(t) $ to an endpoint,
 and then finally comes back to $f(t) $ by passing over
the other endpoint of the given edge.

\indent Now that we have to our disposal the concepts of paths between points, lengths of edges and
internal distances on parts of edges, the {\em length of a path } is
simply defined as the sum of the lengths of the edges and partial edges 
of which the path is made of.  Since the graphs we consider are finite, there are
only finitely many paths between two points $e(s)$ and $f(t)$ and we
can define the {\em  geodesic distance}, $d_{geo}(e(s), f(t))$,
 between these two points as the
minimal length of a path connecting them. If there are no paths
connecting them, then the geodesic distance is defined to be infinite. Each
connected component of the graph is then a compact metric space with
respect to the geodesic distance,  and on each component,
the topology induced by the
metric 
is the same as the
one the component inherits from the ambient compact space.\\
\indent A main emphasis of Kuchement's
research on quantum graphs (see \cite{Ku1, Ku2}) is the study of the spectral properties of a
second order differential operator on a quantum graph. This is aimed at constructing 
mathematical models which can be applied in
chemistry, physics and nanotechnology. On the other hand, our goal is to express geometric features 
of a graph by 
noncommutative geometric means. This requires, in the first place, associating a suitable spectral 
triple to a graph. \\
\indent We have found in the literature several proposals for spectral triples associated to a graph. (See [Co2], 
Section IV.5,  \cite{CMRV, Da, Re}.)  
We agree with Requardt when he states in \cite{Re} that it is a delicate 
matter to call any of the spectral triples the ``right one'', any proposal for a spectral triple associated to a graph 
being in fact determined by the kind of problem one wants to use it for.\\
\indent Our own proposal for a spectral triple is based on the length 
function associated to the edges, which, as was stated before, brings us close to the quantum graph approach. 
A major difference
between the quantum graph approach and ours, however, is that the delicate question of
which boundary conditions one has to impose in order to obtain {\em
  self-adjointness } of the basic differential operator $\frac{1}{i}
\frac{\text{d}}{\text{d}x}$ on each edge disappears in our context. The reason being
that our spectral triple associated with the edges takes care of that issue by
introducing a larger module than just the space of square-integrable functions
on the line.\\
\indent We are now
going to construct a spectral triple for a parameterized graph by
forming a direct sum of all the unbounded Fredholm modules associated to the
edges. This is a special case of the direct sum of curve triples as
studied in Proposition \ref{sumSTr}, so we may state 

\begin{Dfn} \label{dG}
Let  $\cg = ( \cv , \ce) $ be a weighted graph with a parameterization
in a compact metric space $\ct$ and let $\cp$ denote the subset of
$\ct$ consisting of the points of $\cg.$ For
each edge $e = \{e(s)\, |\, 0 \leq s \leq \ell(e)\}$, we let $ST_e$
denote the $e$-triple. The direct sum of the  $ST_e$-triples
over $\ce$ is an unbounded Fredholm module over C$(\ct)$ and a
spectral triple for C$(\cp)$, which we call the {\bf graph triple}
of $\cg$ and denote
$ST_{\cg}$. The metric induced by the graph triple on the set of
points $\cp$ is denoted $d_{\cg}$.
\end{Dfn}
Our next result shows that $d_{\cg}$ coincides with the geodesic metric on the graph $\cg$. 

\begin{Prop} \label{dgeo}
Let  $\cg = ( \cv , \ce) $ be a weighted graph with a 
parameterization. Then, for
any points $e(s), f(t)$ on a pair of edges $e, f,$ we have 
\begin{displaymath}
d_{\cg}(e(s), f(t))\, = \, d_{geo}(e(s), f(t)).
\end{displaymath}
\end{Prop}
\begin{proof}
We will define a  set of functions $\cn$ by 
\begin{displaymath}
\cn \,= \, \{ k: \cp  \to \br \, |\, \exists e(s) \, 
\forall f(t):\,  k(f(t))\,=\, d_{geo}(f(t), e(s))\,\}.
\end{displaymath}
\indent
Locally, the function $f \in \cn$ has slope 1 or minus 1 with respect
to the geodesic distance. This is seen in the following way. Let $ f $
be a function in $\cn$, and suppose you move from 
a point, say $f(t_0),
$ on the edge $f$ to the point $f(t_1)$ on the same edge; 
then the geodesic distance to 
$e(s)$ will usually either increase or decrease by the amount $|t_1 -
t_0|$. This means that {\em most often }
 it is expected that in a
neighborhood around a
point $t_0$ the function $k(f(t))$ is given either as 
$k(f(t)) = k(f(t_0)) + (t-t_0)$ or $k(f(t)) = k(f(t_0)) - (t-t_0)$.
On the interval $(t_0 - \delta, t_0 + \delta ) $, the function $t \mapsto
k(f(t))$ is then differentiable and its derivative 
is either
constantly 
$1$ or $-1$. 

\indent
The reason why this picture is not 
 always true is that for  a point $f(t_0) $ on an edge $f$,  the
geodesic paths from $e(s)$ to the points $f(t_0 - \eps)$ and  $f(t_0 +
\eps)$ may use different sets of edges. The function $k(f(t)) $ is
still continuous at $t_0$, but the derivative will change sign from
$1$ to $-1$, or the other way around.
For each edge $f$, the function $k$ on $f$ is then transported via the
homomorphism attached to the $f$-triple onto the function 
$\phi_f(k)\in \ca_{l(f)}$ and further, via the homomorphism $\Phi_{l(f)}$, onto the 
continuous and even function $\Phi_{l(f)}(\phi_f(k))$ on $[-\ell(f), \ell(f)] $ 
on the positive part of the interval. The function $\Phi_{l(f)}(\phi_f(k))$ is also 
piecewise affine  with slopes in the set $\{-1,1\}$. From 
Proposition \ref{cordomD},
we
then deduce that $D_{\ell(f)}\Phi_{l(f)}(\phi_f(k))$ exists and is essentially
bounded by 1.
Given this fact, it follows that for all functions $k$ in $\cn$ and any edge
$f$ we have $\|\,[D_{\ell(f)}, \pi_f(k)]\,\| \leq 1$, so with the
above choice of function $k$, we obtain

\begin{align*} 
& \forall e,f \in \ce, \, \forall t \in [0, \ell(f)], \, \forall s \in [0,
\ell(e)]: \\ & d_{geo}(f(t), e(s) ) \,= \,k(f(t)) - k(e(s))\, \leq
\, d_{\cg}(f(t), e(s)).
\end{align*} 
\indent
Suppose now that we are given points as before, $e(s)$ and $f(t)$, on 
some  edges $e$ and $f$. Then there exists a geodesic path from $e(s)$
to $f(t)$, and we may without loss of generality assume that this path
is an ordinary graph-theoretic path consisting of the edges $e_1,
\dots , e_k$ such that the starting point of $e_1$ is $e(s)$ and the
endpoint of $e_k$ is $f(t)$. Any continuous function $g$ on $\cp$
which has the property that 
\begin{displaymath}
\underset{e \in \ce}{\max}\,\|\,[D_{\ell(e)}, \pi_e(g)]\,\|\, \leq \,
1
\end{displaymath} 
also has the property that its derivative on each edge 
 is essentially  
bounded by 1. In particular, this means that for each edge
$e_j,$ with $ 1 \leq j \leq k$, we must have $|g(e_j(\ell(e_j)))- g(e_j(0))| \leq
\ell(e_j) $ and then, since $e_i(\ell(e_i)) = e_{i+1}(0) $ for $1 \leq i
\leq k-1$, we have

\begin{align*}
|g(e(s))-g( f(t))| =& |g(e_1(0))- g(e_k(\ell(e_k)))|\\ 
 \leq & \sum_{j=1}^k |g(e_j(0)) - g(e_j(\ell(e_j)))|\\
 \leq & \sum_{j=1}^k \ell(e_j)\\
 = & d_{geo}(e(s), f(t)).
\end{align*}
\indent
Since this holds for any such function $g$, we deduce that 
\begin{displaymath}
d_{\cg}(e(s),f(t)) \,  \leq \, d_{geo}(e(s),f(t)),\end{displaymath}
 and the proposition
follows.
  \end{proof}

\smallskip
\indent
It should be remarked that the unbounded Fredholm modules of 
noncommutative geometry can be refined to give more information about the
topological structure of the graph. Suppose for instance that the
graph is connected but is not a tree. It then contains at least one
closed path, {\em a cycle}. For each cycle, one can obtain a
parameterization by adding the given parameterizations or the inverses
of the given parameterizations of the edges which go into the
cycle. In this way, the cycle can be described via a spectral triple
for a circle of length equal to the sum of the edges used in the
cycle. If this unbounded Fredholm module is added to the direct sum of 
the curve-triples coming
from each edge, then the 
 geodesic distance is still measured by the new
spectral triple, and this triple will induce an element in the K-homology of the 
graph, as in  \cite{Co2}, Chapter IV, Section $8. \delta$.
This element in the K-homology group
will be able to measure the winding number of a
nonzero continuous function around this cycle. One may, of 
course, take
one such a summand for each cycle and in this way obtain an unbounded Fredholm module
which keeps track of the connectedness type of the graph. We
have not yet found a suitable use for this observation in the case of finite
graphs, but we would like to pursue this idea 
in connection with our later
study of the Sierpinski gasket in Section 8. (See, especially, item 
(vi) in Section 9, along with the relevant discussion following it.)

\section{Infinite trees}

\noindent
In the first place, it is not possible to create unbounded Fredholm modules for
infinite graphs by taking a direct sum of the unbounded Fredholm modules
corresponding to each of the edges and each of the cycles. The problem
is that there may be an infinite number of edges of length bigger than
some $\delta > 0$. In such a case, the direct sum of the Dirac
operators will not have a compact resolvent any more. We will avoid
this problem by only considering graphs which we call {\em finitely
  summable trees} and which we define below, but before that, we want to remark
that the concept of a path remains the same, namely, a finite collection of
edges leading from one vertex to another without repetitions of
vertices.  
\begin{Dfn}
An infinite  graph $\cg = (\cv, \ce ) $ is a {\bf finitely summable tree} 
 if the following
conditions are satisfied:
\begin{itemize}
\item[(i)] There are at most countably many edges.
\item[(ii)] There exists a length function $\ell$ on $\ce$ such that
  for any edge $e$, we have $\ell(e) > 0. $
\item[(iii)] For any two vertices $u, v $ from $\cv$, there exists
  exactly one path between them.  
\item[(iv)] There exists a real number $p \geq 1$ such that $\sum_{e
    \in \ce} \ell(e)^p < \infty$. (In that case, the tree
is said to be {\bf {\em p}-summable}.) 
\end{itemize}
\end{Dfn}

\indent An infinite tree which is not finitely summable may create
problems when one  tries to look at it in the same  way as we did for finite
graphs. It may, for instance, 
not be possible to  embed such a graph into a locally compact space in a
reasonable way. To indicate what the problems may be, we now discuss an
example.  Think 
 of the bounded tree whose vertices $v_n$  are indexed by
$\bn_0$ and whose edges all have length 1 and are given by  $e_n
:= \{ \, \{v_0, v_n \} \, |\, n \in \bn\}, $  i.e. all the edges go
out from $v_0$ and the other vertices are all endpoints, and they are
all one unit of length away. If the
points $v_n$ for $n \geq 1 $ 
have to be distributed in a symmetric way in a metric
space, then the distance between any two of them should be the same and 
no subsequence of this sequence will be convergent. This
shows right away that the graph can not be embedded in a compact space in a
reasonable way, but it also shows that the point $v_0$ does not have a
compact neighborhood, so even an embedding in a locally compact space
 is not possible. The restriction of working with finitely summable trees will be
 shown to be sufficient in order to embed the graph in a compact
 metric space. Before we show this, we would like to mention that the
 point of view of this article is to consider edges as line segments
 rather than as pairs of vertices. 
We will start, however, from the
 graph-theoretic concept of edges as pairs of vertices  and then show very
 concretely that we  can obtain a model of the corresponding
 simplicial complex, even in such a way that the edges are continuous
 curves in a compact metric space.

 \begin{Dfn} \label{parinfgraph} A finitely summable tree $\cg$ with
   vertices $\cv$ and edges $\ce$ is said to have a {\bf parameterization}
   if  $\cv$ can be represented injectively as a subset of a metric
   space $(\ct, d)$ and for each edge $e = \{p, q\}$, there exists a
   continuous curve, $\{e(t) \in \ct \, |\, 0 \leq t \leq \ell(e)\}$,
   without self-intersections and leading from one of the vertices of $e$
   to the other. Two curves  $e(t)$ and $f(s)$  for different 
 edges $e$ and $f$ may only intersect at endpoints. 
The length of the defining interval, $\ell(e)$, for a curve
   $e(t),$ is called the {\bf length of the curve}. The set of points $e(t)$ on
   all the curves is denoted $\cp$. 
\end{Dfn}

\indent
For an infinite  tree, there is only one path between any
two vertices  in $\cv.$ This is the basis for the following concrete
parameterization of an infinite $p$-summable tree inside a compact
subset of  the Banach space $\ell^p(\ce)$. We first fix a vertex $u$ in
$\cv$ and map this to $0 $ in $\ell^p(\ce)$; then the unique path
from $ u $ to any other vertex $v$ determines the embedding in a 
canonical way.  In order to describe
this construction, we will denote by $\delta_e$ the canonical unit basis vector in $\ell^p(\ce)$
corresponding to an edge $e$ in $\ce$. 

\begin{Prop}
Let $\cg = (\cv, \ce)$ be a finitely summable tree, $u$ a vertex in
$\cv$ and $T_u: \cv \to \ell^p(\ce)$ be defined by 
\begin{displaymath}
T_u(w) : = \begin{cases} 0\\ 
\underset{j=1}{\overset{k}{\sum}} \ell(e_j)\delta_{e_j}
\end{cases}
\quad
\begin{array}{l}
\text{ if }\, w = u ,\\
 \\  \text{ if } \,
w \neq u  \,\text{ and  
  the path  is } \, \{e_1, \dots , e_k\}.
\end{array}
\end{displaymath}
For an arbitrary edge $e =\{w_1, w_2\}$, we choose an orientation such that the
first vertex is the one which is nearest to $u$.    
Let us suppose that $w_1$ is the first vertex.  Then, given such an oriented edge  
$e =(w_1,w_2)$, a curve $e_u: [0, \ell(e)] \to \ell^p(\ce)$ is defined by 
\begin{displaymath}
e_u(t) \, : =  \, T_u(w_1) \, + \, t\delta_e.
\end{displaymath}
\indent
Let $\cp_u$ denote the set of points on all the curves $e_u(t)$ and
let $\ct_u$ denote the closure of $\cp_u.$ Then $\ct_u$ is a norm
compact subset of $\ell^p(\ce)$ and    
the triple  $(\ct_u, \,  T_u, \, \{e_u(t) \, | \, e \in
\ce\})$ constitutes a parameterization of $\cg$, in the sense of Definition 7.2.
\end{Prop} 

\begin{proof}
The only thing which really has to be proved is that the set $\ct_u$
is a norm compact subset of $\ell^p(\ce).$ Let then $\eps > 0 $ be
given and choose a finite and connected subgraph  $\cg_0 = (\cv_0,
\ce_0) $ of $\cg$ such that $\sum_{e \in \ce_0}\ell(e)^p + \eps^p \, >
\, \sum_{e \in \ce}\ell(e)^p $ and $u$ is a vertex in $\cv_0$.  
We denote by $\cp_{u,\ce_0}$ the set of points on all the curves 
$e_u(t)$, where $e\in \ce_0$. We then remark that $\cp_{u,\ce_0}$ is just a finite union of
compact sets, and is therefore compact. By construction, any point in $\cp_u$ is
within a distance $\eps$ of $\cp_{u,\ce_0}$, so the closure $\ct_u$ of $\cp_u$ is
compact.
\end{proof}

\indent
It is not difficult to see that for two different vertices $u$
and $v$, the metric spaces $\ct_u$ and $\ct_v$ obtained via
the above construction 
are actually isometrically isomorphic
via the rather natural identification of the sets $\cp_u$ and $\cp_v$ described just below.
\begin{Prop}\label{tutv}
Let $\cg = (\cv, \ce)$ be a finitely summable tree. Consider $u$ and $v$ two different vertices in
$\cv$ and let $T_u: \cv \to \ell^p(\ce)$, respectively $T_v: \cv \to \ell^p(\ce)$, be defined as in Proposition 7.3. 
Then the  mapping $S_{uv}: \cp_u \rightarrow \cp_v$ defined by
\begin{displaymath}
S_{uv}(x) : = \begin{cases} T_v\circ T_u^{-1}(x) \quad \mathrm{ if }\, x\in T_u(\cv),\\
T_v(w_j)+t \delta_e  \quad \mathrm{ if } \,
x =e_u(t), \, e=(w_1,w_2),\, t\in (0,l(e)),
\end{cases} 
\end{displaymath}
where $w_j$ is the nearest vertex to $v$ 
among $w_1$ and $w_2$,
is an isometric isomorphism of $\cp_u$ onto $\cp_v$.  
\end{Prop} 
\begin{proof}
 Given two points, say $x, y$, in $\cp_u$, we have to show that
 \begin{displaymath}
 d(S_{uv}(x), S_{uv}(y)) \, = \, d(x,y).
\end{displaymath}
 There has to be two different
 arguments showing this, according to the cases where $x$ and $y$ are
 on the same or on different edges. We will only consider the case
 where $x = e_u(s)$ on an
 edge $e$  and $y = g_u(t)$ on a different  edge $g.$ 
Then there exists a finite
 set of vertices $w_1, \dots,  w_k$ such that the closed line segments in $\ell^p(\ce)$ given by 
\begin{displaymath}
[e_u(s), T_u(w_1)], \, [T_u(w_1), \, T_u(w_2)], ... ,  [T_u(w_{k-1}),
\, T_u(w_k)], \,  [T_u(w_k), \, g_u(t)], 
\end{displaymath}
all belong to $\cp_u$ and constitute the unique path herein from $x$ to
$y.$ The distance in $\cp_u$ from $x$ to $y$ is given by 
\begin{align*}
&\|x-y\|_p \,= \\
& \left( \| T_u(w_1)- e_u(s)\|^p + \| g_u(t) -
  T_u(w_k)\|^p + \sum_{i=1}^{k-1} \ell(\{w_i, w_{i+1}\})^p  \right) 
^{1/p}.  
\end{align*}
The term $\| T_u(w_1)- e_u(s)\|$ is either  $s$ or $\ell(e) - s$. It
is $s$ if $w_1$ is closer to $u$ than $x$, and $\ell(e) - s$ otherwise.
The distance in $\cp_v$ from $S_{uv}(x)$ to $S_{uv}(y)$ is given by
\begin{align*}
&\| S_{uv}(x)-S_{uv}(y) \|_p \,= \\
& \left( \| T_v(w_1)- e_v(s')\|^p + \| g_v(t') -
  T_v(w_k) \|^p + \sum_{i=1}^{k-1} \ell(\{w_i, w_{i+1}\})^p  \right) 
^{1/p}.  
\end{align*}
A moment's reflection will make it clear to the reader
that the distance between $S_{uv}(x)$ and $S_{uv}(y)$ is the same
as the distance between $x$ and $y.$ Thus $S_{uv} $  defines an
isometry between the sets of points $\cp_u$ and $\cp_v$.  
\end{proof}
Since the sets $\cp_u$ and $\cp_v$
are dense in $\ct_u$ and $\ct_v$, respectively, we see that $S_{uv} $
extends to a natural isometry between $\ct_u$ and $\ct_v $. The compact
metric spaces $\ct_u$ and $\ct_v$ are then isometrically isomorphic
via an isometry which commutes with the parameterizations, so we may
introduce 
\begin{Dfn}
Let $\cg \, = \, ( \cv, \ce)$ be a $p$-summable infinite tree, $u$ 
a vertex in $\cv$, and  $(\ct_u, \,  T_u, \, \{e_u(t) \, | \, e \in
\ce\})$ the parameterization of $\cg$ introduced in Proposition 
7.3. In view of Proposition 7.4, we may use the simpler notations $\ct,$ 
respectively $\cp,$ 
to denote the sets $\ct_u$ and $\cp_u$.  
(In the sequel, we will also use the notation $\ct^p, \cp^p $ when needed.)
The parameterization 
is called the $\mathbf{p}$-{\bf parameterization} of $\cg$. 
The metric given by $\Vert \cdot \Vert_p  $ is denoted 
by $d_p$ and the boundary of $\cp$ in $\ct$ is denoted by $\cb$.  
\end{Dfn}

\smallskip
\indent It is clear that if an
infinite tree is $p$-summable, then it is also $q$-summable
for any real number $q>p$.
A natural question is, of course, if
the $p$- and the $q$-parameterizations are homeomorphic or even
Lipschitz equivalent as metric spaces. We can easily show that the
spaces are homeomorphic but the Lipschitz equivalence may not be
automatic, although it is easy to establish in many concrete
examples. We have found a condition which ensures Lipschitz
equivalence; it is nearly a
tautology, but rather handy for the study of concrete examples. To be able to express
this property, we first observe that the $p$-parameterization has a
concrete realization via a base vertex $u$ as $\ct_u$. This is a
subset of $\ell^p(\ce) $ and since the $p$-norm here dominates the
$\infty $-norm, we can use this norm to introduce a metric $d_{\infty} $
on $\ct_u$ by $d_{\infty}(x,y)\,  := \, \|x-y\|_{\infty}.$ 
If one goes back to the
description of $\ct_u$, one can realize that $d_{\infty}(x,y) $ is
independent of $u$ and, in fact, depends only on the unique path
(may be even infinite in both directions)
which leads from $x$ to $y.$

\begin{Prop} \label{homeo}
Let $\cg = (\cv, \ce) $ be a $p$-summable infinite tree,
and $q$ be a real number such that $q > p $. Further, let $\ct^p, \ct^q$ denote the $p$-, respectively, 
$q$-parameterization compact metric spaces of $\cg$.  Then these spaces are
 homeomorphic via the natural embedding of $\ct_p$ into $\ct_q.$

\indent Moreover, the metric spaces $\ct^p$ and $\ct^q$ are Lipschitz equivalent if there
exists a constant $k > 0 $ such that
\begin{displaymath} 
\forall  x, y \in \ct^p: d_p(x,y) \, \leq \, k d_{\infty}(x,y).
\end{displaymath}
\end{Prop} 

\begin{proof}
Let us take a base vertex $u$ and consider the concrete representations
of $\ct^p$ and $\ct^q$ as $\ct_u^p$ and $\ct_u^g$. Since $p < q$, we
have $\ell^p(\ce) \subset \ell^q(\ce).$ Let $\iota$ denote the
canonical embedding; then $\iota$ is a contraction and consequently,
$\iota(\ct_u^p) $ is a compact subset of $\ct_u^q$ containing the
point set $\cp_u^q$ of $\ct_u^q$. Hence, $\iota$ induces a
homeomorphism between the two compact spaces $\ct_u^p $ and
$\ct_u^q.$
 \smallskip
\noindent
Let us now assume that  
\begin{displaymath} 
\forall  x, y \in \ct_u^p: d_p(x,y) \, \leq \, k d_{\infty}(x,y),
\end{displaymath}
and continue to work inside $\ct_u^p.$ Then we get
\begin{displaymath} 
\forall  x, y \in \ct_u^p: d_q(x,y) \leq  d_p(x,y)
\leq  k d_{\infty}(x,y) \leq k d_q(x,y),
\end{displaymath}
so the metrics are Lipschitz equivalent.
\end{proof}
\smallskip
\indent
It took us quite some time to realize how different various parameterizations
may be. Later in this section, the reader can see (in Figure 3) a picture of the fractal
which is usually associated to the free noncommutative group on 2
generators. In Connes' book \cite{Co2}, on page 341, one can see quite a different
picture. In the first case, the boundary is totally disconnected,
whereas in the Poincare disk picture used in \cite{Co2}, the boundary is
the unit circle. The $p$-parameterization has the property that it
separates different boundary points very much.

\begin{Prop} \label{boundary}
Let $\cg \, = (\cv, \, \ce)$ be an infinite $p$-summable tree and
$\ct$ its $p$-parameterization space. The set $\cb$ of boundary points
of $\cp$ in $\ct$ is closed and totally disconnected.
\end{Prop}
\begin{proof}
Let $a$ and $b$ be two different  points in $\cb$ and let $\delta$ be
a positive number less than $d(a,b)/3$. To $\delta$  we associate a
finite connected  subgraph $\cg_0 \, = \, (\cv_0, \, \ce_0 )$ such that 
\begin{displaymath} 
\sum_{e \in \ce_0} \ell(e)^p + \delta^p \, > \, \sum_{e \in \ce}
\ell(e)^p.
\end{displaymath}
The graph $\cg_0$ is a finite tree, so it has at least two ends,
i.e. vertices of degree 1. By examination of the concrete space
$\ct_u$, one can easily see that the set of points  $\cp_0$ of $\cg_0$ 
 has nonempty interior
as a subset of $\ct$, and that this interior, say $\overset{\circ}{\cp_0}$,
is exactly all of $\cp_0$ except its endpoints. Given a point $x$ in
$\cp$  which
is in the complement, say $\cc$, of $\overset{\circ}{\cp_0}$ in $\ct$, 
we see that there must be a path from $x$ to $\cp_0$ which
meets $\cp_0$ at an endpoint. This means that all the points in
$\cc$ are grouped into a finite number of
pairwise disjoint sets according to which end vertex in $\cv_0 $ is
the nearest. Suppose now that $x$ and $y$ are points in $\cc$ such
that their nearest end vertices in $\cv_0,$ say $v$ and $w$, respectively,
are different. Then, by the
construction of the metric space $(\ct, d)$, one finds that $d(x,y) >
d(v, w) > 0.$ This shows that the connected components of $\cc$  can
be labeled by the end vertices of $\cv_0$ as $\cc_v$
 and there is a positive distance
between any two different components. 

\smallskip
\indent Let us  now return to the boundary points $a, b$ and show that they
fall in  different components  $\cc_u$ and $\cc_v$. Suppose to the contrary
that both $a$ and $b$ belong to the same
component, say $\cc_u$. Then there exist vertices $v$  and $w$ in
$\cc_u \cap \cv $ such that $d(a,v) < \delta$ and $d(b,w) < \delta. $ Since $v$
and $w$ are in the same component $\cc_u$, the unique path between them
must be entirely in $\cc_u$. This means that it uses none of the edges
from $\ce_0$; so, by the construction of $\cg_0$,  we get $d(u,v) <
\delta $ and then $d(a,b) < 3\delta < d(a,b),$ a contradiction. Since the components induce a
covering of the boundary by open sets, we deduce that
$a$ and $b$ are in different components of the boundary and thus that the boundary
is totally disconnected. The closedness of $\cb$  follows, as we can take an
increasing sequence of finite connected subgraphs of $\cg,$ say
$(\cg_n)$, such that the union of the sequence of open sets
$(\overset{\circ}{\cp_n})$ in  $\ct$ is all of $\cp.$  
\end{proof}

\smallskip
\indent We are not going to study the properties of the $p$-parameterization in
more details since our main  interest is to describe certain aspects of
graphs with the help of noncommuative geometry. The possibility of
embedding a graph as a dense subset of a compact metric space $\ct$
shows that we can study this space through a spectral triple which is
the sum of triples associated the edges of the tree. The ambient metric space
has a metric which is constructed such that the space becomes
compact. On the other hand the natural distance between points on a
tree is the length of the unique path between the points and this will
quite often not be a bounded metric on the tree, so we will refer to
this distance as the geodesic distance.

\begin{Dfn}
Let  $\cg = (\cv, \, \ce) $ be a graph which is embedded in a metric
space $\ct$ such that each edge, $ e = (x_0,
x_1),$ where   $x_0$ and $x_1$ are vertices in $\cv,$ has a
parameterization $r_e(t)$ with $r_e(0) = x_0$ and $r_e(\ell(e) )=
x_1$. Then the geodesic distance between any two points $x, y  \in
\ct$, which are  situated on edges on the graph, is defined 
to be the infimum of the sums of the lengths of the  corresponding
intervals which have to be used in order to construct a curve  from
$x$ to $y$ based on the parameterization given for each edge.
\end{Dfn}

The $p$-parameterization
provides a framework which makes it possible to obtain a spectral triple as
a sum of triples based on curves, where the curves are parameterized
edges of the tree. We will concentrate our investigations on the
spectral triple which can now be obtained in the same way as was
done for finite graphs. Before we embark on this, we would like to
mention that unless the summability number $p$ equals 1, one should
not expect that the geodesic distance on $\cp$ is bounded. Think, for instance,
of a tree embedded in $\br_+$ with vertices $(v_n)$ indexed by the
natural numbers, edges of
the form $\{v_n, v_{n+1}\}$ and lengths $\ell(\{v_n, v_{n+1}\}) = 1/n.$
This graph is 2-summable and the set of path lengths is unbounded.
Let us further remark that, by applying the triangle inequality, one can see
that the geodesic distance is always larger than
the distance on the $p$-parameterization 
space, which is a subset of $\ell^p(\ce)$. Thus it may very well
be that the geodesic distance is unbounded on $\cp$ and hence, in such a case,
 is not extendable to a continuous function on the compact space $\ct.$

\indent
In \cite{Co2}, Section IV.5, {\em Fredholm Modules and Rank-One Discrete
  Groups}, Connes studies modules which are associated  to trees. His
modules are $\ell^2$-spaces over sets consisting of points and edges,
whereas the ones we are going to construct are sums of $L^2$-spaces
over the edges, as described for finite graphs in Section 6. There may be closer
relations than we can see right now between the two types of modules;
at least, it seems that Proposition 6 on page 344 of \cite{Co2} is
related to our Theorem \ref{infsptrip}. 

\begin{Dfn} \label{ratrip}
Let $\ct$ be a compact space,   $r: [0, \alpha] \to \ct $ a
continuous curve, and let  $a$ be a real number. The $r,a$-triple
\begin{displaymath}
(\mathrm{C}(\ct), H_{\alpha}, D_\alpha^a ) 
\end{displaymath}
is defined as
a translated unbounded Fredholm module (in the sense of Definition
\ref{transpt})  of the $r$-triple defined in  
Definition \ref{DefSTr}.   The operator $D_\alpha ^a$ is the translate of 
$D_\alpha$ given by 
 $D_{\alpha} + aI$.
\end{Dfn}

\indent Since for any $f \in \mathrm{C}(\ct)$ we have $[D_{\alpha} + aI,
\pi_{\alpha}(f) ] \, = \, [D_{\alpha},
\pi_{\alpha}(f) ]$, the change in the Dirac operator does not affect
the spectral triple much, but the eigenvalues of the Dirac
operator all  get translated by the number $a$.

\begin{Thm} \label{infsptrip}
Let  $\cg = (\cv, \ce)$ be
 a $p$-summable infinite tree with $p$-parameterization in the compact metric
 space $\ct.$  The direct sum over all the
edges $e$ from $\ce$ of the
unbounded Fredholm modules
\begin{displaymath}  
 (\mathrm{C}(\ct), H_{\ell(e)},
 D_{\ell(e)}^{\pi/(2\ell(e))})
 \end{displaymath}
is a spectral triple for $\mathrm{C}(\ct)$, which is denoted 
 $ST_{\cg}: = ( \mathrm{C}(\ct), H_{\cg}, D_{\cg}).$ 
It can only be a
finitely summable module for a real number $s > 1$. Further,
for a given  $s > 1$, it is finitely summable 
 if and only if 
\begin{displaymath}
\sum_{e \in \ce} \ell(e)^s \, < \, \infty.
\end{displaymath}
Moreover, the metric $d_{\cg}$ induced by $ST_{\cg}$  on the points
$\cp$ of the infinite
tree $\cg$ is the geodesic distance. 
\end{Thm} 
\begin{proof}
The proof relies in many ways on the corresponding proof for finite graphs. 
For points from the set $\cp, $ we shall see that 
the previous arguments can be reused to a large extent.  The real
problem occurs when we have to deal with the boundary points $\cb$ in $
\ct,$ i.e. the set of points 
$\ct$ which are not in $\cp.$ We noticed above
that the geodesic distance may not be extended to a continuous function
on $\ct,$ so we can not just copy the proof of Proposition \ref{dgeo} and
define the set $\cn$ similarly. 
Instead, we consider the set, say $\cf_\cg$, 
 of finite connected subgraphs $\cg_0
=(\cv_0, \ce_0)$  of $\cg.$ In the proof of Proposition \ref{boundary},
we saw that the complement of the open set $\overset{\circ}{\cp_0}$
in  $\ct$ consists of a finite collection of pairwise disjoint closed sets $\cc_v $
labeled by the endpoints of $\cg_0.$ 
This makes it
possible to define a dense algebra of continuous functions on $\ct$
which will have uniformly bounded commutators for all of  our
unbounded Fredholm modules associated with the edges. We simply define
the  set of
functions $\cn$ by looking at functions which have bounded commutators
for the unbounded Fredholm module of Definition \ref{dG} applied to some
$\cg_0$ in $\cf_\cg$  and
are constant on each of the components, outside
$\overset{\circ}{\cp_0}$,
 i.e. the value of
such a function $f$ on a component $\cc_v$ is $f(v).$
Any function $f$ in $\cn$ must have a bounded commutator with $D_{\cg}$, since the commutator is zero except at edges in a
subgraph $\cg_0,$ and here it is supposed to have bounded commutator
with the Dirac operator associated to $\cg_0$. The functions in $\cn$
constitute a self-adjoint algebra of continuous functions on $\ct$
 since the set $\cf_\cg$ is upwards directed under
inclusion. Ordinary points in $\cp$ can be separated by functions from
$\cn$ in the same way as this was done in the proof of Proposition
\ref{sumSTr}. For a boundary point $b$ and an ordinary point $p $ in $\cp$,
it will always be possible to find a graph $\cg_0 $ in $\cf_\cg$ such
that $p$ is an inner  point in the set of points $\cp_0$ associated to
$\cg_0$. The function $k$, which on $\cp_0$ is defined by $k(x) :=
d_{geo}(x, p)$ and continued by constancy on the components associated to
the end vertices of $\cg_0$, belongs to $\cn $ and separates $p$ and the
boundary point $b.$ The reason being that $p$ is assumed to be an
inner point in $\cp_0$, so it will have positive geodesic distance to
any  end vertex 
of $\cg_0.$ For two different boundary points $b$ and $c$, one
may use the proof of Proposition \ref{boundary} to obtain a subgraph
$\cg_0$ in $\cf_\cg$ such that $b$ and $c$ fall in different
components. Suppose  $b $ is in $\cc_v$  and $c $ is in $\cc_w$ for end vertices 
$v$ and $w$ of $\cg_0$. As seen above, there is a function $f$
in $\cn$ such that $f(v) \neq f(w)$; but $f(b) = f(v) $ and $f(c) =
f(w)$, so $f$ separates $b$ and $c$. 

In order to prove that we have a spectral triple  for C$(\ct)$, 
we then have to show that the resolvents
of $D_\cg$ are compact, when bounded. This will be proven below, but first
we will show that the metric induced by the set $\{f \in
\mathrm{C}(\ct)  \, | \, \|[D_{\cg}, \pi(f)]\| \leq 1 \, \}$ is the
geodesic distance on $\cp.$ We will then return to the
analogous problem for finite graphs in the proof of  Proposition
\ref{dgeo}. Again, we show that for any two points $a$ and $b$ from
$\cp$, we may find a connected subgraph $\cg_0$ in $\cf_{\cg}$  such that its set of
points in $\ct$ contain both $a$ and $b$. We then conclude as in the {\em
  finite }  case that $d_{\cg}(a,b)  \, = \,  d_{geo}(a,b). $\\
\smallskip
\indent In order to see that we have an unbounded Fredholm module and prove the
summability statement, we have to look at the eigenvalues of
\begin{displaymath}
\underset{e \in \ce}{\oplus} D_{\ell(e)}^{\pi/(2\ell(e))}.
\end{displaymath}
The eigenvalues for each summand form the set $\{(k+1/2)\pi/\ell(e)\,
|\, k \in \bz \, \}$; so $(D_{\cg}^2 +I)^{-1} $ has the  following
doubly indexed set of eigenvalues 
\begin{displaymath}
 \quad \bigg(\frac{4\ell(e)^2}{(2k+1)^2\pi^2 +4\ell(e)^2}\,\bigg)_{(
   \, k \in \bz, \, e \in \ce \,) }.
\end{displaymath}
\indent
Remark that the  same eigenvalue may occur, for different edges,  but only 
 a finite number of times. 
Since 
\begin{displaymath}
 \frac{4\ell(e)^2}{(2k+1)^2\pi^2 +4\ell(e)^2} \leq
\frac{4\ell(e)^2}{(2k+1)^2}
\end{displaymath}
 and $\sum_{e \in \ce} \ell(e )^s <  \infty $,  
 we see that $(D_{\cg}^2+I)^{-1}$ is compact.

\indent
With respect to the summability, we consider  for a real number $s>0$ the sum
\begin{displaymath}
\sum_{e \in \ce} \sum_{k \in \bz} \bigg(\frac{4\ell(e)^2}{(2k+1)^2\pi^2
  +4\ell(e)^2}\bigg)^{s/2}.
\end{displaymath} 
The term $(2k+1)^2$ in the denominator implies that 
 we must have $s >1 $ in order to obtain a finite sum.
Now, for $s>1$, we rewrite the double sum as follows:
\begin{align*}
& 2^s\cdot \sum_{e \in \ce} \ell(e)^s \cdot \sum_{k \in \bz} \left (\frac{1}{(2k+1)^2\pi^2
  +4\ell(e)^2}\right)^{s/2}=\\
& 2^{s+1}\cdot \sum_{e \in \ce} \ell(e)^s \cdot 
\sum_{k \in \bn_0} \left (\frac{1}{(2k+1)^2\pi^2
  +4\ell(e)^2} \right)^{s/2}.
\end{align*}
For each $e \in \ce$, there is a $k_e \in \bn_0 $ such that $(2k+1)\pi > 2\ell(e)$ for every 
$k > k_e$. Then, for each $e\in \ce$, we obtain that 
\begin{align*}  
& \sum_{k=0}^{k_e}\left ( \frac{1}{(2k+1)^2\pi^2
  +4\ell(e)^2}\right)^{s/2} +
 \left ( \frac{1}{2}\right )^{s/2}\sum_{k=k_e+1}^{\infty}\left ( \frac{1}{(2k+1)^2\pi^2
  }\right)^{s/2}
 \leq \\ 
& \sum_{k\in \bn_0} \left ( \frac{1}{(2k+1)^2\pi^2
  +4\ell(e)^2}\right)^{s/2}  \leq  
\sum_{k\in \bn_0} \left ( \frac{1}{(2k+1)\pi
  }\right)^s   
\end{align*}  
and hence it follows that the module is $s$-summable if and only if $s>1$
and the
sum $\sum_{e \in \ce} \ell(e)^s$ is finite.  
\end{proof}

\begin{exmp}  \label{Caylf2}
The Cayley graph for the noncommutative free group on
2 generators $\baf_2$.
\end{exmp} 
\indent
There is a nice description of this graph as a
fractal. We start with the neutral element $e$ at
the origin of $\br^2$. Then the generators $\{a,b, a^{-1},
b^{-1}\}$ are placed on the axes at the points
\begin{displaymath} 
\{ (1/2, 0), (0, 1/2),
(-1/2,0), (0, -1/2)\}.
\end{displaymath}

\begin{figure}[h]
  \psfrag{a}{$\scriptstyle a$}
  \psfrag{-a}{$\scriptstyle a^{-1}$}
  \psfrag{b}{$\scriptstyle b$}
  \psfrag{-b}{$\scriptstyle b^{-1}$}
  \psfrag{e}{$\scriptstyle e$}
  \psfrag{ab}{$\scriptstyle ab$}
  \psfrag{a-b}{$\scriptstyle ab^{-1}$}
  \psfrag{a^2}{$\scriptstyle a^2$}
  \psfrag{ab-a}{$\scriptstyle aba^{-1}$}
  \psfrag{aba}{$\scriptstyle aba$}
  \psfrag{ab^2}{$\scriptstyle ab^2$}
  \psfrag{a^2b}{$\scriptstyle a^2b$}
  \psfrag{a^2-b}{$\scriptstyle a^2b^{-1}$}
  \psfrag{a^3}{$\scriptstyle a^3$}
  \psfrag{a-b-a}{$\scriptstyle ab^{-1}a^{-1}$}
  \psfrag{a-ba}{$\scriptstyle ab^{-1}a$}
  \psfrag{a-b^2}{$\scriptstyle ab^{-2}$}
 \includegraphics[scale=0.4]{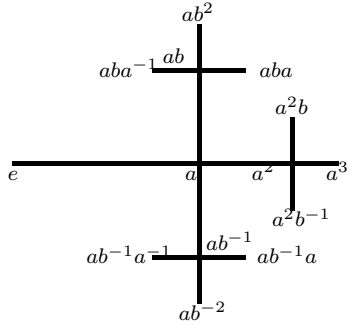}
   \caption{A portion of the Cayley graph of $\baf_2$.}\label{Cayley1}
  \end{figure}

Traveling right along an edge represents multiplying on the right by $a$, while
traveling up corresponds to multiplying by $b$. Each new edge is drawn
at half size of the previous one to give a  fractal image. 
We start by illustrating  in Figure \ref{Cayley1} how all words,
say $CGa$,  
which begin with an $a$ are positioned on the tree.

\begin{figure}[h]
 \includegraphics[scale=0.30]{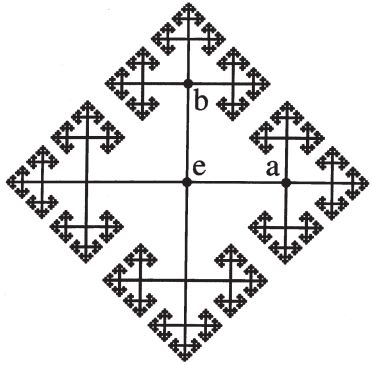}
 \caption{The Cayley graph of $\baf_2$, viewed as a fractal tree.}\label{Cayley2}
  \end{figure}

\smallskip
\indent Figure \ref{Cayley2} shows the  entire graph which consists of
four identical fractal images. 
Each one, say  $CG_a, CG_{a^{-1}}, CG_b, CG_{b^{-1}}$, represents all
the words  starting 
with $a$, $a^{-1}$, $b$ or  $b^{-1}$.

This forms an infinite tree
with 4 edges of length 1/2, 12 edges of length
1/4, 36 edges of length 1/8 and, generally, $4\cdot 3^{n-1}$ edges of length
$2^{-n}$. The sum $\sum_{e\in \ce} \ell(e)^s$ from  Theorem \ref{infsptrip} can then be written as 
\begin{displaymath}
\sum_{e \in \ce} \ell(e)^s \, = \,
\sum_{m=1}^{\infty}4\cdot 3^{m-1}2^{-ms}\,
=\,\frac{4}{3}\sum_{m=1}^{\infty}(3\cdot 2^{-s})^{m}.
\end{displaymath}
So the module is finitely summable if and only if $s >\log 3  / \log 2$.
Hence, the spectral triple of $\baf_2$ has metric dimension $\log 3  / \log 2$,
and this is exactly the Minkowski and also the Hausdorff dimension of
the closure of 
the Cayley graph of $\baf_2$. 

\smallskip
\indent  As was seen just above, given any 
real number $p > \log 3 / \log 2$, 
the graph is $p$-summable; so we can consider its
$p$-parameterization $\ct^p.$ Since the lengths of the edges decrease
geometrically like $2^{-n}$, it follows that  we have
\begin{displaymath}
\forall x, y \in \ct^p: \quad d_p(x,y) \leq 2 d_{\infty}(x,y).
\end{displaymath} 
Hence, by Proposition \ref{homeo},  $\ct^q$ and $\ct^p$ are Lipschitz
equivalent for any $q > p > \frac{\log 3}{\log 2}.$ Further, some computations
with the closure of the Cayley graph for $\baf_2$ constructed above
show  
that the metric coming from $\br^2$ on this set is also Lipschitz
equivalent to 
the $d_{\infty}$ metric. Therefore, the $\br^2$ fractal is actually Lipschitz
equivalent to any of the $\ct^p$ parameterization spaces for $p >
\frac{\log 3}{\log 2}.$

\smallskip 
\indent On page 341 of \cite{Co2}, one can find a representation
of the Cayley graph of $\baf_2$ as a subset of the Poincare disc, and
as we mentioned above, the boundaries of these two different
parameterizations of the Cayley graph of $\baf_2$ are very different.

\subsection{Complex Dimensions of Trees}

The mathematical theory of complex fractal dimensions finds its origins
in the study of the geometry and spectra of fractal drums \cite{La1, La3} and, in particular,
of fractal strings \cite{LM,LP}. In the latter case, it is developed in the research monograph
\cite{L-vF1} and significantly further expanded in the recent book \cite{L-vF2}. We establish here
some connections between this theory and our work. 
\smallskip

\noindent \begin{Dfn}\label{ZetaF}
Let  $p$ be a real number greater than $1$ and $\cg$ an infinite {\em p}-summable tree  with vertices $\cv$
and a countable number of parameterized edges $\ce = \{e_n \, |\, n \in
\bn \}$.  Assume that the lengths of the edges $\ell(e_n) $ converges to $0$ as
$n \to \infty$. The
{\bf zeta function} of the tree $\cg$ is denoted $\zeta_\cg(z)$ and is defined for $\rm{Re}(z)>p$ by
\begin{displaymath}
\zeta_\cg(z)={\rm tr}\left ( \vert D_\cg \vert ^{-z}\right ).
\end{displaymath} 
\end{Dfn}

In view of Theorem 7.10, $\vert D_\cg \vert$ has the following doubly indexed set of eigenvalues.
\begin{displaymath}
\left( \frac{(2k+1)\pi}{2\ell(e)}\right)_{(
   \, k \in \bn_0, \, e \in \ce \,) },
\end{displaymath}
each with multiplicity 2. \\
\indent Let $\zeta(z)$ denote the Riemann zeta function. By writing
\begin{align*} 
\zeta(z) \,= &\,\sum_{l=1}^{\infty} l^{-z} = \sum_{k=0}^{\infty}(2k+1)^{-z}
+ \sum_{k=1}^{\infty} (2k)^{-z}\\
 =&\,  \sum_{k=0}^{\infty}(2k+1)^{-z}
+ 2^{-z} \cdot \zeta(z),
\end{align*} 
we obtain    
\begin{equation} 
\sum_{k=0}^{\infty}(2k+1)^{-z} =  (1-2^{-z})\cdot \zeta(z).
\end{equation}  
\indent  Hence, for $\rm{Re}(z)>p$, we deduce that 
\begin{align*} 
{\rm tr}(|D_{\cg}|^{-z} ) \, =& \, 2\sum_{e \in \ce} \sum_{k \in \bn_0
} \left( \frac{2\ell(e)}{(2k+1)\pi}\right)^z\\
=& \, \frac{2^{z+1}}{\pi^z} \cdot \left ( \sum_{e \in \ce}\ell(e)^z \right) \cdot \left( \sum_{k=0}^{\infty}(2k+1)^{-z}\right ) \\
=& \, \frac{2^{z+1}}{\pi^z} \cdot (1-2^{-z})\cdot \left ( \sum_{e \in \ce}\ell(e)^z \right ) \cdot \zeta(z).
\end{align*}

\begin{exmp}\label{ZetaCayley}
Let $\cg$ be the Cayley graph of $\baf_2$. Then, the zeta function of $\cg$ 
has a meromorphic extension to all of $\bc$ given by 
\begin{displaymath}
\zeta_\cg(z)=\frac{8}{\pi^z}\frac{1-2^{-z}}{1-3\cdot 2^{-z}}\cdot \zeta(z), \text{ for z } \in \bc.
\end{displaymath} 
Indeed, this is true for ${\rm Re}(z)>\log 3 /\log 2$, by the last displayed
equation in Example 7.11. Hence, by analytic continuation, it is true for all
$z\in \bc$. \\
\indent Aside from a trivial multiplicative factor $f(z)$, which is an entire function, the zeta function of $\cg$
is of the same form as the spectral zeta function of a self-similar fractal string,
which is always equal to the product of the Riemann zeta function and the geometric zeta function
of the fractal string (see \cite {La2, La3, LM, LP} and Chapter
1 in \cite{L-vF1} or \cite{L-vF2}).
\end{exmp}
\indent We refer to \cite{L-vF2}, Chapters 1, 4 and 5, for the precise definition 
of the complex dimensions of a fractal string, using the notions of ``screen and window".

\begin{Dfn}\label{ComplDim}
Assume that $\zeta_\cg$ admits a meromorphic continuation to an open neighborhood
of a window $W \subset \bc$. Then the {\bf visible complex dimensions} of $\cg$
(relative to $W$) are the poles in $W$ of the meromorphic continuation of $\,\zeta_\cg$.
The resulting set of visible complex dimension is denoted by $\frak{D}_\cg(W)$: 
\begin{displaymath}
\frak{D}_\cg(W):=\{z \in W\mid \zeta_\cg \text{ has a pole at } z \}.
\end{displaymath} 
If $W=\bc,$ then we simply write $\frak{D}_\cg$ for $\frak{D}_\cg(\bc)$
and call the elements of $\frak{D}_\cg$ the {\bf complex dimensions} of $\, \cg$. 
\end{Dfn}
\indent It follows that the set of complex dimensions of $\cg$, the Cayley graph of $\baf_2$, 
is given by 
\begin{displaymath}
\frak{D}_\cg= \{1\}\cup \{ \cd_\cg +\sqrt{-1}\cdot k\cdot {\bf p}\mid k \in \bz\},
\end{displaymath}
where, in the terminology of \cite{L-vF1, L-vF2}, $\cd_\cg := \log 3/\log 2 $ 
is the Minkowski dimension of $\cg$ and ${\bf p}:=2\pi/\log 2$ is its oscillatory period. \\
\indent Note that this is analogous both to what happens for self-similar strings (see \cite{L-vF2}, 
Section 2.3.1 and especially Chapter 3, where such strings are allowed to have dimension greater than 1),
as was mentioned earlier, and for the Sierpinski drum (see \cite{L-vF2}, Section 6.6.1), viewed as a 
fractal spray (more specifically, as the bounded, infinitely connected planar domain with boundary 
the Sierpinski gasket). Indeed, it follows from (\cite{L-vF2}, Eqs. (6.81) and (6.82)) that the set of (spectral) complex dimensions
of the Sierpinski drum is equal to $\{2\} \cup \frak{D}_\cg$, with $\frak{D}_\cg$ as in the last displayed equation. In our present 
situation, the value 1 appears naturally since it corresponds to the dimension
of any edge of the tree, while in \cite{L-vF2}, the additional value 2 occurs because the Sierpinski 
drum is viewed as embedded in $\mathbb{R}^2$.\\

\begin{rem}
More generally, choose a suitable window $W$ contained in the half-plane ${\rm Re}(z)>0$. We deduce 
from the discussion following Definition \ref{ZetaF} that the zeta function 
$\zeta_\cg(z)$ of any $p$-summable tree $\cg$ (as in Theorem \ref{infsptrip}) is given by
\begin{displaymath}
\zeta_\cg(z):=g(z)\cdot \zeta_\cl(z)\cdot \zeta(z)=g(z)\cdot\zeta_\nu(z).
\end{displaymath}
Here, $g(z)$ is an entire function which is nowhere vanishing in $W$ and $\zeta_\cl(z)$ is the
meromorphic continuation to $W$ of the geometric zeta function $\sum_{e \in \ce}\ell(e)^z$
of the fractal string $\cl=\cl_\cg:=\{{\ell_e}\}_{e \in \ce}$ associated with $\cg$.
Furthermore, $\zeta_\nu(z):=\zeta_\cl(z)\cdot \zeta(z)$ is the spectral zeta function
of $\cl$ (\cite{La3,LP}, and \cite{L-vF2}, Theorem 1.19).
In particular, we have  
\begin{displaymath}
\frak{D}_\cg(W)=\{ 1\} \cup \frak{D}_\cl(W),
\end{displaymath}
where $\frak{D}_\cl(W)$ is the set of visible complex dimension of $\cl$ (as in Definition \ref{ComplDim}). 
\end{rem}

\bigskip
\section{The Sierpinski Gasket: Hausdorff Measure and Geodesic Metric} \label{sgsection}
\begin{figure}[h]
  \includegraphics[scale=0.2]{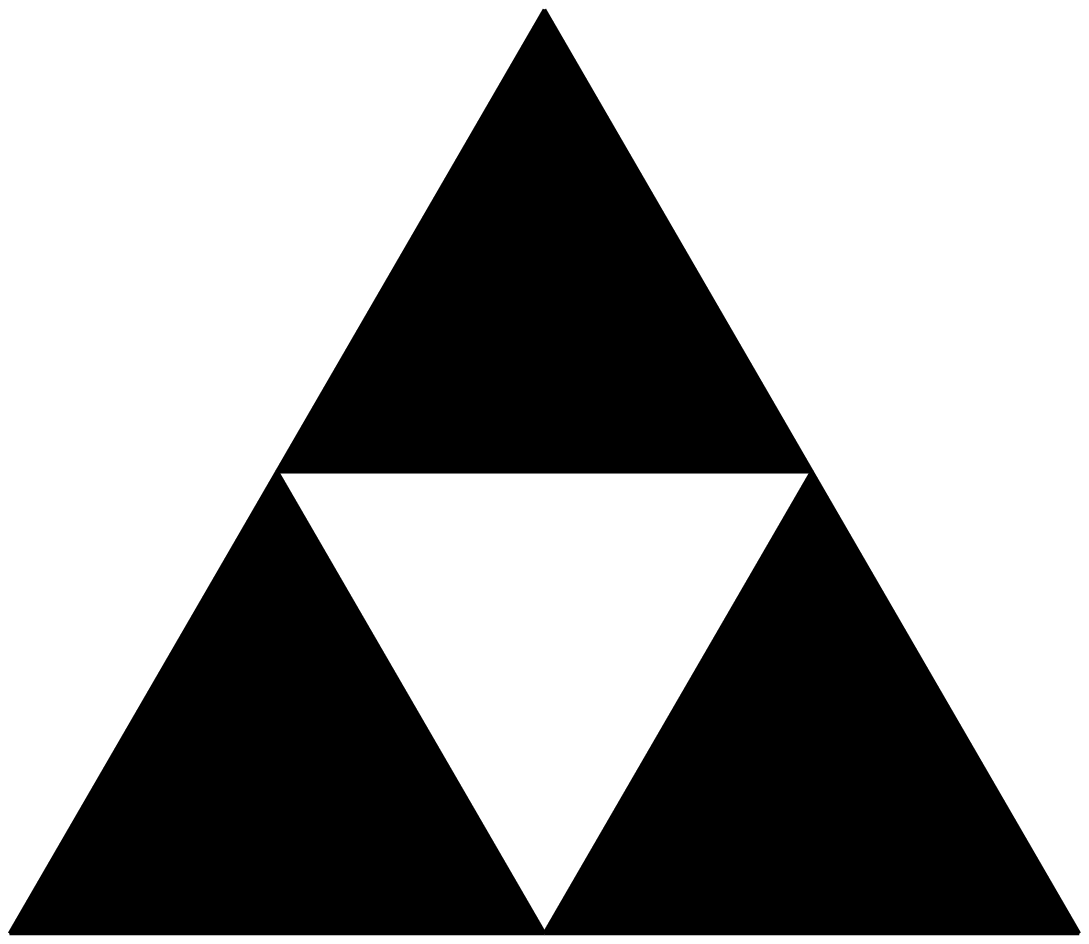}
  \includegraphics[scale=0.2]{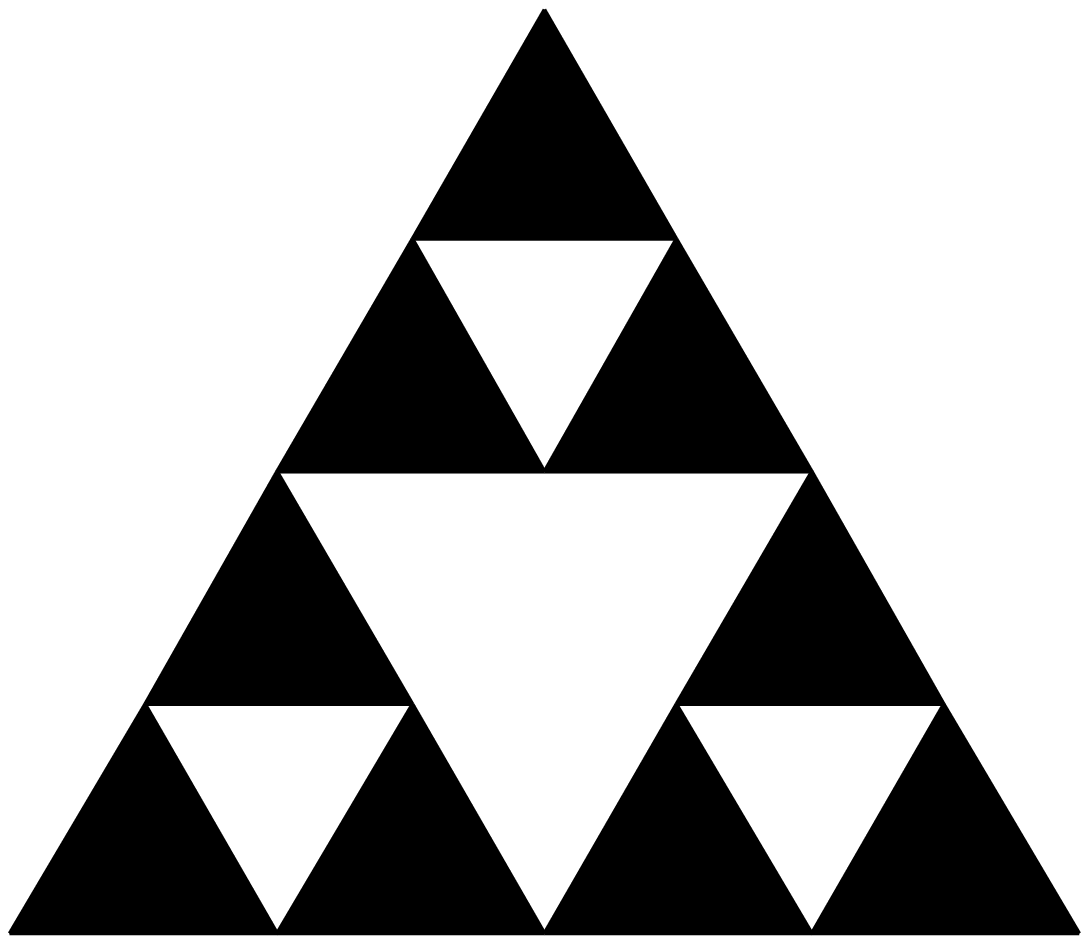}
  \includegraphics[scale=0.2]{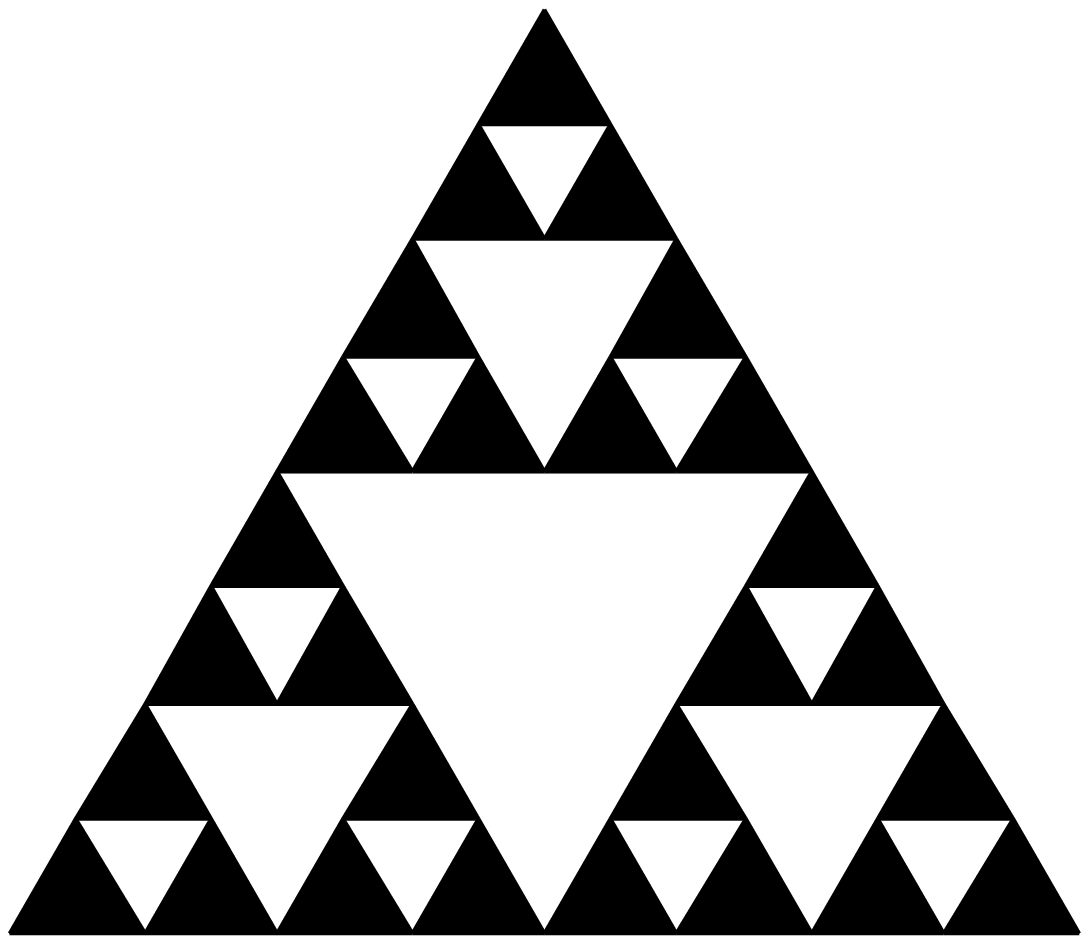}
  \caption{The first 3 steps in the standard construction  of the Sierpinski gasket.}\label{SG2}
  \end{figure}

\noindent The Sierpinski gasket is well known and is described in many places. 
In particular, it is a connected fractal subset of the Euclidean
plane $\br^2$; in fact, it can be viewed as a continuous
planar curve which is nowhere differentiable.  
We refer the reader to the books by Barlow, Edgar, and  Falconer \cite{Ba,Edg1, 
Edg2, Fa}, which all contain good descriptions and a lot of information about
this fractal set.
The gasket  can be obtained in many ways. The most common is probably
the one, illustrated in Figure \ref{SG2}, 
 where one starts with a solid equilateral triangle in the plane
 and cut out
 one open equilateral triangle of half size, and then continue to cut
 out open triangles of smaller and smaller sizes. Hence, in the $n\text{-th}$ step, one 
cuts away $3^{n-1} $
open equilateral triangles with side length equal to  $2^{-n}$  that of the
original one. 

\indent In our present study, we would rather like to consider a
construction where 
the Sierpinski gasket, denoted from now on by $\cs\cg$, is obtained as
the closure of the limit of an 
increasing sequence of 
sets in
$\br^2$. This
means that we take as a starting point an equilateral triangle, but
not solid anymore; just its border, consisting of the 3 sides, and the
three vertices. We call this figure $\cs\cg_0$. 
The next figure, $\cs\cg_1$, is obtained by adding another triangle of half size,
and turned upside down relative to  $\cs\cg_0$, and so on. This procedure is
well known, and illustrated in  Figure
\ref{SG1} in the introduction.  We are not so much interested in the algorithm 
used to construct the Sierpinski gasket. Instead, our goal is to describe the
topological properties of the gasket via noncommutative methods. From
this point of view, it seems better to describe 
$\cs\cg_1$ as consisting of  3 equilateral triangles of half size {\em and
glued together at 3 points.} Clearly, the set $\cs\cg_0$ is still in a natural
way  a subset of $\cs\cg_1$. The following figure $\cs\cg_2$ then consists of
$3^2$ triangles of size $2^{-2}$ of $\cs\cg_0$, and they are glued
together at $3 + 3^2$ points. And finally, $\cs\cg_n$
consists of $3^n$ triangles, each shrinked to a size $2^{-n}$ of the
starting one, and glued together at $3(3^n-1)/2$ points.   The figures
\ref{SG1} and \ref{SG2}
illustrate, when compared,  the well known fact that the closure of the union of the
sets $\cs\cg_n$ equals the intersection of the decreasing sequence of
the sets, which are obtained via the cutting procedure.

There are, of course, identifications between points in
the figures such that each figure $\cs\cg_n$ is a subset of the next one 
$\cs\cg_{n+1}$, but
this will be taken care of by the construction of a spectral
triple for C$(\cs\cg)$ as a direct sum of triples associated to
each of the triangles which appear in any of the figures $\cs\cg_n$.
The advantage of this way of constructing a spectral triple for the
Sierpinski gasket is
that it keeps track of the {\em holes} in the gasket. 

\indent We will use  a circle as the basis for the construction of
an unbounded Fredholm module for each of the triangles contained in  $\cs\cg_n$, and
then ultimately 
obtain a spectral triple for the Sierpinski gasket. In order to avoid using 
too many $\pi$'s, we fix the side length of $\cs\cg_0$ to $2 \pi / 3$. 
The perimeter of each triangle which is used to form $\cs\cg_n$ is then $2
\pi / 2^n$. We will think of each of the triangles shown above as having
a horizontal edge. Further, we will introduce a natural
parameterization of each triangle in $\cs\cg_n$ by using the right-hand
corner on the bottom line as the starting point and then run counterclockwise
by arclength. In this way, we get for each such  triangle in $\cs\cg_n$
an isometry of the circle of 
radius $2^{-n}$  onto this triangle, when both are 
equipped with arclength as metric.   
We will now introduce a numbering of the triangles which go into this
construction and also define the associated spectral triples. 
\begin{Dfn} \label{Deltan}
\begin{itemize}
\item[(i)] Given $n $ in $\bn_0$, choose a numbering of the $3^n$
  triangles of size $2^{-n}$  which form $\cs\cg_n$, and let
  $\Delta_{n,i}$, $i \in \{1, 2, \dots , 3^n\}$, denote the
  numbered triangles.
\item[(ii)] Let, for each $n$ in $\bn_0$ and each $i$ in $\{1, \dots,
  3^n\}$, the mapping $r_{n,i} : [-2^{-n}\pi, 2^{-n}\pi ] \to
  \Delta_{n,i}$ be defined such that $r_{n,i}(0) $ equals the lower
  right-hand corner of $\Delta_{n,i}$ and the mapping is an isometry,
  modulo $2^{1-n}\pi$, of this interval onto the triangle $\Delta_{n,i}
  $ equip-ped with the geodesic distance as metric, and the counterclockwise orientation. 
The mapping $r_{n,i} $ induces a surjective
  homo-morphism $\Phi_{n,i} $ of C$(\cs\cg)$ onto 
C$([-2^{-n}\pi, 2^{-n}\pi ])$ by
\begin{displaymath}
\forall t \in [-2^{-n}\pi, 2^{-n}\pi ],\, \, \forall f \in
\mathrm{C}(\cs\cg): \quad \Phi_{n,i}(f)(t) \, :=\, f(r_{n,i}(t)). 
\end{displaymath}
\item[(iii)] Let, for each $n$ in $\bn_0$ and each $i$ in $\{1, \dots,
 3^n\}$, the unbounded Fredholm module $ST_{n,i}(\cs\cg)\, := \, (\mathrm{C}(\cs\cg), \, H_{n,i}, \, D_{n,i})
$ for $\cs\cg$ be given by
\begin{itemize}
\item[(1)] $H_{n,i} \,:= \, H_{2^{-n}\pi} \quad \text{(see Definition 3.1)}$;
\item[(2)] the representation $\pi_{n,i}:$ C$(\cs\cg)\to
  B(H_{n,i})$ is defined for $f$ in C$(\cs\cg)$ as the
  multiplication operator which multiplies by the function
  $\Phi_{n,i}(f)$.
 \item[(3)] $D_{n,i} \, :=\, D^t_{2^{-n}\pi}\quad \text{(see just 
 after Definition 2.1)}$.
 \end{itemize}
\end{itemize}
\end{Dfn}

\begin{Thm} \label{SgST}
The direct sum of all the unbounded Fredholm modules  $ST_{n,i}(\cs\cg)$ for $n$
in $\bn_0$ and $i$ in $1, \dots , 3^n$ gives a spectral triple for 
$\cs \cg$. This spectral triple is denoted $ST(\cs\cg) =
(\mathrm{C}(\cs \cg), H_{\cs \cg}, D_{\cs \cg}) $,  
and it  is $s$-summable if and only if $s > \log 3 / \log 2$. Hence, its 
metric dimension is equal to $\log 3 / \log 2$, which is also the
Minkowski (as well as the Hausdorff) dimension of the Sierpinski gasket.
\end{Thm}

\begin{proof}
As usual, we have to prove that the algebra  of continuous functions on
$\cs\cg$, which have bounded commutators with the  Dirac
operator $D_{\cs\cg}$, is dense in C$(\cs\cg)$. As before, it is enough to show that this
algebra separates the points in $\cs\cg$. In this case, our task is
easy. We just remark that the real-valued linear functions on $\br^2$ separate
points and we will prove that any real-valued linear functional on $\br^2$, say
$\psi(x,y) = ax+ by$, has the property that its restriction to $\cs\cg$
has a bounded commutator with the sum of the Dirac operators. 
We will think of each of the triangles shown above as having an horizontal edge, as before. 
Moreover, we will consider a Euclidean coordinate system such that the $x$-axis is also horizontal. 
Let us
then consider the composition of $\psi $ with 
any of the parameterization mappings $r_{n,i}$. This 
yields a continuous $2^{1-n}\pi$-periodic function, say $f$, on the interval
$[-2^{-n}\pi, 2^{-n}\pi ]$ or rather on $\br$, 
which is affine on the 
intervals $[-2^{1-n} \pi /3, 0], [0, 2^{1-n} \pi /3 ], [-2^{-n} \pi, -2^{1-n} \pi /3]$, and 
$[2^{1-n} \pi /3, 2^{-n}\pi ]$.
The slopes of all the functions $\psi\circ r_{n,i} $ will belong to the set 
\begin{displaymath}
S:=\left \{ a, \frac{-a+b\sqrt{3}}{2}, \frac{-a-b\sqrt{3}}{2} \right \}.
\end{displaymath}
Note that $S$ is independent of both $n$ and $i$.  
According to Lemma 2.3, for each pair $n,i$ the  function $\psi \circ
r_{n,i}$ 
is in the domain of $D^t_{n,i}$ and $\|[\pi_{n,i}(f), D^t_{n,i}]\|
\leq \max(S).$ Hence, $\|[\pi_{\cs\cg }(f), D_{\cs\cg }]\|
\leq \max(S)$ and we have proven that the restriction
 of any affine
function f on $\br^2$ to $\cs\cg$ has the property that
$\|[\pi_{\cs\cg}(f), D_{\cs\cg}]\|$ is bounded. 

\smallskip
\indent In order to check that  $(D_{\cs\cg}^2+I)^{-1}$ is compact and establish
its summability properties, we  compute
for each $n \in \bn_0$ the set of  eigenvalues $E_n$  of
$D_{\cs\cg}$, which are added on  in the $n$-th step. 
\begin{itemize}
\item[$n=0:$] $E_0 = \{ (2k+1)/2 \,|\,k \in \bz\}$, each of multiplicity
  $1 = 3^0.$
\item[$n > 0:$] $E_n = \{  2^{n-1}(2k+1) \,|\,k \in \bz\}$, each of
  multiplicity  $3^n.$
\end{itemize}

\noindent (Note that the formula for $E_n (n>0)$ is valid for $n=0$ as well. However, by scaling,
it is naturally derived from that for $E_0$.)

\indent Now that we know the eigenvalues, we can compute the trace
  of the operator $|D_{\cs\cg}|^{-z}$ for a given complex
  number $z$, and we will in the next theorem show that the trace  is finite
  for   Re$(z) > \log 3 / \log 2 $.
\end{proof}

\begin{rem}
It is well known that the Sierpinski gasket has Minkowski and Hausdorff dimensions
equal to $\log 3/ \log 2$. Indeed, it follows from the fact that $\cs\cg$ is a self-similar
set satisfying the Open Set Condition and can be constructed out of 
3 similarity transformations each of scaling ratio $1/2$. (See e.g. \cite{Fa}, Chapter 9.)
\end{rem}

\begin{Thm} \label{Tr}
  Let $\zeta(z)$ denote the Riemann zeta function. Then, for any complex
  number $z$ such that $ \mathrm{Re}(z) > \frac{\log 3}{\log 2}$, 
\begin{displaymath}
{\rm tr}\left ( \vert D_{\cs \cg}\vert ^{-z}\right ) \, = \, 2^{z+1}\cdot\frac{1-2^{-z}
  }{1 - 3\cdot 2^{-z}}\cdot \zeta(z).
\end{displaymath}
Moreover, for any Dixmier trace $\mathrm{Tr}_{\omega}$ on $B(H_{\cs\cg})$, we have  
\begin{displaymath} 
\mathrm{Tr}_{\omega}\left (\vert D_{\cs \cg}\vert ^{- \frac{\log 3}{\log
    2}} \right ) \, = \, \frac{4}{\log 3}\cdot \zeta \left (\frac{\log 3}{\log2 } \right ).
\end{displaymath}
\end{Thm}
\begin{proof}
We recall here the formula 1 from page 32: 
\begin{displaymath} 
\sum_{k=0}^{\infty}(2k+1)^{-z} =  (1-2^{-z})\cdot \zeta(z). 
\end{displaymath} 
Hence, for Re$(z)>\log 3/ \log 2$ we deduce that
\begin{align*} 
{\rm tr}(|D_{\cs \cg}|^{-z} ) \, =& \, \sum_{n \in \bn_0} \sum_{k \in \bz
}3^n 2^{-nz}|k+1/2|^{-z} \\
=& \, \sum_{n \in \bn_0} 3^n 2^{-nz} \sum_{k \in \bz}2^z|2k+1|^{-z}\\
=& \, \frac{2^z}{1 - 3 \cdot 2^{-z}}\cdot   2 \cdot \sum_{k=0}^{\infty} |2k+1|^{-z} \\
=& \, 2^{z+1}\cdot \frac{ 1-2^{-z}}{1 - 3\cdot 2^{-z}} \cdot \zeta(z).
\end{align*}
Then, by Proposition 4 on page 306 in \cite{Co2}, we obtain 
\begin{align*} 
 \mathrm{Tr}_{\omega}\left (\vert D_{\cs \cg}\vert ^{- \frac{\log 3}{\log
    2}}\right  )\,  = \,& \underset{x \to 1+}{\lim} (x-1){\rm tr}\left ( \left (
  \vert D_{\cs \cg}\vert ^{-\frac{\log 3}{\log 2}}\right )^x \right) \\ = \,&
 \underset{x \to 1+}{\lim} (x-1)\cdot 2^{x\frac{\log3}{\log 2}+1}\cdot \frac{1-2^{-x\frac{\log3}{\log 2}}
 }{1 - 3\cdot 2^{-x\frac{\log3}{\log 2}}}\cdot \zeta \left (x\frac{\log3}{\log
  2}\right )\\
= \,&
 \underset{x \to 1+}{\lim} \frac{x-1}{1- 3^{1-x}}2(3^x -1)\cdot 
\zeta\left (x\frac{\log3}{\log
  2}\right )\\
= \,& \frac{4}{\log3}\cdot \zeta\left (\frac{\log3}{\log 2}\right ). 
\end{align*}

We note that it is also possible to compute the classical limit
expression for the Dixmier trace of a measurable operator, obtained as
the limit as $N\rightarrow \infty$ of $(1/\log
N)\sum_{j=0}^{N-1}|\lambda_j|^{-(\log 3 / \log 2)}$, where
$(\lambda_j)_{j=1}^{\infty}$ are the characteristic values (here, the
eigenvalues) of the given compact operator, written in non increasing
order.
  \end{proof}

\begin{rem}
The zeta function of the Sierpinski gasket can now be defined for $
{\rm Re}(z) > \frac{\log 3 }{\log2} $  
by $\zeta_{SG}(z):={\rm tr}\left ( \vert D_{SG}\vert ^{-z}\right ).$
  In view of Theorem \ref{Tr}, it
 has a meromorphic continuation to all of $\mathbb{C}$ and is given by 
\begin{displaymath}
\zeta_{SG}(z)=
2^{z+1}\cdot \frac{1-2^{-z} }{1-3\cdot 2^{-z}}\cdot \zeta(z), \text{ for } z\in \bc.
\end{displaymath}
It also follows from the foregoing formula that 
the zeta function of the Sierpinski gasket (as defined  
just above) and the zeta function for the Cayley 
tree given in Example \ref{ZetaCayley} are proportional modulo the function
$h(z)=\pi^z\cdot 2^{z-2}$. Hence, they have the same set of complex dimensions (see Section 7.1), given by 
\begin{displaymath}
\frak{D}_{SG}= \{1\}\cup \{ \cd_{SG} +\sqrt{-1}\cdot k\cdot {\bf p}\mid k \in \bz\},
\end{displaymath}  
where $\cd_{SG} := \log 3/\log 2 $ 
is the Minkowski dimension of the gasket and ${\bf p}:=2\pi/\log 2$ is its oscillatory period. 
This is in agreement with a conjecture made in \cite{La4}, Section 8, when discussing 
the `geometric complex dimensions' of the gasket. We note, however, that the value 1 was not included 
in \cite{La4} but appears naturally in our context since it corresponds to the dimension
of the boundary of any of the holes (circles or triangles) of the gasket. Following \cite{La4},
we may refer to $\frak{D}_{SG}$ as the set of {\em geometric complex dimensions} of the Sierpinski gasket.    
\end{rem}

\bigskip
\indent A Dixmier trace, $\mathrm{Tr}_{\omega}$,  on $B(H_{\cs\cg})$ induces a
positive linear functional,  $\tau$, on C$(\cs\cg)$, as
stated in \cite{Co2}, Proposition 5, Chapter IV.2, and proved in \cite{CiGS}. 
It turns out that  $\tau $ is a nonzero multiple of the
integral with respect to the Hausdorff measure, say $\mu$,  on the
Sierpinski gasket (see e.g. \cite{Ki}, \cite{KiL1}, or \cite{St}), 
which, in turn, has a very natural
description in terms of functional analytic concepts. 
To show this, we will first give a description of the Hausdorff integral,
which can be used to establish the proportionality between $\tau$ 
and the Hausdorff measure (or rather, integral). 
Having this to our disposal, we can deduce from the second part of 
Theorem \ref{Tr}  that for any
continuous function $f$ on the gasket, we have
\begin{displaymath}
\tau(f)\, := \,
\mathrm{Tr}_{\omega}\left ( \pi_{\cs\cg}(f)\vert D_{\cs \cg}\vert ^{-
  \frac{\log 3}{\log  2}} \right )\, = \, 
\frac{4}{\log 3}\cdot \zeta\left (\frac{\log 3}{\log2
}\right ) \cdot \int_{\cs\cg}f(x)\mathrm{d}\mu(x).
\end{displaymath}
In order to establish this relation, we will use the description of the gasket as an 
 increasing sequence
$\cs\cg_n$ of graphs, but this time we will, for a given nonnegative 
integer $n$, focus on
the set consisting of the $3^{n+1} $ midpoints of the sides in the
$3^n$ triangles $\Delta_{n,i}$ of size $2^{1-n}\pi /3 $ contained in $\cs\cg_n.$ 
We refer the reader to Figure 1 in the introduction for the pictures of the first few sets $\cs\cg_n $ and
now begin the formal description.\\
\smallskip 
\indent The triangles in $\cs\cg_n$ are denoted
$\{\Delta_{n,i}\, | \, 1 \leq i \leq 3^n \}$ and we denote  
the midpoints of a 
triangle $\Delta_{n,i} $ by  $\{x_{n,i, j} \, |\,
1 \leq j \leq 3 \}. $ For any natural number $n$, we then define a positive
linear functional $\psi_n$  of norm 1, a state, on C$(\cs\cg)$  by
\begin{displaymath} 
\forall \, f \, \in \mathrm{C}(\cs\cg): \quad \psi_n(f)\, :=\, 
3^{-(n+1)}\sum_{i=1}^{3^n}\sum_{j=1}^3 f(x_{n,i,j}).  
\end{displaymath}

\begin{Prop} \label{psi}
Let $\mu$ denote the Hausdorff probability measure on the Sierpinski
gasket $\cs\cg$ and let $\psi$ denote the state on
$\mathrm{C}(\cs\cg)$ defined by 
\begin{displaymath}
\forall \, f \, \in \mathrm{C}(\cs\cg): \quad \psi(f)\,:=\, \int_{\cs\cg}f(t)\mathrm{d}\mu(t).
\end{displaymath}
   Then the sequence of states $(\psi_n)_{n \in \bn} $ converges to
   $\psi $ in the weak*-topology on the dual of $\mathrm{C}(\cs\cg).$
\end{Prop}

\begin{proof} 
We will first show that for each complex continuous function $f $
on $\cs\cg$, 
the sequence $(\psi_n(f))_n$ is  a Cauchy sequence. This follows 
from the fact that any such function $f$ is uniformly
continuous and the points $x_{n,i,j}$ are evenly distributed on the
gasket. To be more precise, let $\eps > 0 $ be given, then 
  choose, by the uniform
continuity of $f$, an $n_0$ in
$ \bn$ so large that for each $h$ in $\{1, \dots , 3^{n_0}\}$ and any two points
$x, y$ in $\cs\cg$ which are inside or on the triangle 
 $\Delta_{n_0,h } $, we have $|f(x)-f(y)| \leq \eps.$ Let then $n$ be a
natural number bigger than $n_0$ and  let us  consider  the average, say
$\nu^n_{n_0,h}(f),$ of the values
$f(x_{n,i,j}) $ over all the midpoints $x_{n,i,j}$ from  
 $\cs\cg_n$ which are on the
border or inside the triangle $\Delta_{n_0,h}$. Then the estimate for the 
 variation
of $f$ over the points inside or on $\Delta_{n_0,h} $ implies that 
\begin{displaymath}
\left \vert \nu^n_{n_0,h}(f) - \frac{f(x_{n_0,h,1}) + f(x_{n_0,h,2}) +f(x_{n_0,h,3})}{3} \right \vert
      \, \leq \, \eps.
\end{displaymath}

\indent We can then establish the Cauchy property by using the following inequalities:

\begin{align*} 
\forall n \geq n_0: \, & | \psi_n(f) - \psi_{n_0}(f) | \\
 = \, \, & \left \vert  3^{-n_0} \sum_{h = 1}^{3^{n_0}} \left (\nu_{n_0,h}^n(f) - \frac{f(x_{n_0,h,1}) +
   f(x_{n_0,h,2}) +f(x_{n_0,h,3})}{3}\right )\right \vert \\
 \leq & \,  \eps.
\end{align*} 

\indent This shows that the sequence $(\psi_n)_{n \in \bn}$ is weak*-convergent 
to a state, say $\phi,$ on C$(\cs\cg).$ By
construction, it is clear that $\phi$ is self-similar with respect to
the basic affine contractions which define Sierpinski's gasket, as
described for instance by Strichartz in \cite{St}. Hence, $\phi = \psi$
and the proposition follows.
 \end{proof}
Having this to our disposal, we can establish the claimed 
relationship between $\tau$ and $\psi.$  
\begin{Thm} \label{proportional}
Let $\tau$ be the functional on C$(\cs\cg)$ given by  
\begin{displaymath} \tau(f)\, : = \,
\mathrm{Tr}_{\omega}\left ( \pi_{\cs\cg}(f)\vert D_{\cs \cg}\vert
^{-\frac{\log 3}{\log  2}} \right ) 
\end{displaymath} 
and $\mu$ the Hausdorff
probability measure on $\cs\cg$. Then, for any continuous complex-valued
function $f$ on $\cs\cg$, we have
\begin{displaymath}
\tau(f)\, = \, 
\frac{4}{\log 3}\cdot \zeta\left ( \frac{\log 3}{\log2
}\right ) \cdot \int_{\cs\cg}f(x)\mathrm{d}\mu(x) \, = \, \frac{4}{\log 3}\cdot \zeta \left ( \frac{\log 3}{\log2
}\right ) \cdot \psi(f).
\end{displaymath}
\end{Thm}
\begin{proof}
Let $f$ be  a continuous {\em real-valued} function on $\cs\cg$ and
$\eps >0 $ a positive real number. Let us then go back to the proof of Proposition
\ref{psi} and choose $n_0 \in \bn$ and, for $n > n_0$ and $h \in \{1,
\dots , 3^{n_0} \}$, define $\nu^n_{n_0,h}(f)$ as above.
We then restrict our attention to the portion of the Sierpinski gasket which is contained inside or on the 
triangle $\Delta_{n_0, h}$ and denote this space by $\cs\cg_{n_0,h}$. It
follows that for the identity function on $\cs\cg_{n_0,h}$, say $I_{n_0,h} $,
and for $f_{n_0, h}$ the analogous restriction of $f$, we have in the
natural ordering on C$(\cs\cg_{n_0,h})$,
\begin{displaymath}
(\nu^n_{n_0,h}(f) - \eps) I_{n_0,h} \, \leq \,  f_{n_0,h} \, \leq \,
(\nu^n_{n_0,h}(f) + \eps) I_{n_0,h}.
\end{displaymath}
For each $h \in  \{1, \dots , 3^{n_0} \}$, we can naturally define
a spectral triple for $\cs\cg_{n_0,h}$ by deleting all
the summands of $ST(\cs\cg)$ which are based on triangles outside
$\Delta_{n_0,h}.$ To any such triple, we can associate a corresponding
functional $\tau_{n_0,h}$ and we get 
\begin{displaymath}
\tau(f) \, = \, \sum_{h=1}^{3^{n_0}}\tau_{n_0,h}(f_{n_0,h}) \quad
\text{and} \quad  \tau_{n_0,h}(I_{n_0,h})\, = \, 3^{-n_0}\tau(I). 
\end{displaymath}
Since $\tau$  is a positive functional, the inequalities above give
\begin{displaymath} 
\sum_{h=1}^{3^{n_0}} (\nu^n_{n_0,h}(f) - \eps)(3^{-n_0}\tau(I)) \, 
\leq \, \tau(f)\, \leq \,  \sum_{h=1}^{3^{n_0}} (\nu^n_{n_0,h}(f) +
\eps)(3^{-n_0}\tau(I)),
\end{displaymath}  
where $I$ is the identity for the algebra of continuous functions on $\cs\cg.$ 
If we then go back to the proof of Proposition  \ref{psi}, we see that
\begin{displaymath} 
| \psi(f) - 3^{-n_0}\sum_{h=1}^{3^{n_0}} \nu^n_{n_0,h}(f) | \, \leq \, \eps,
\text{ so } |\tau(I)\psi(f) - \tau(f)| \, \leq \, \tau(I)\eps.
\end{displaymath}  By the second part of Theorem
\ref{Tr}, we know that 
\begin{displaymath}
\tau(I) \, = \, \frac{4}{\log
  3}\cdot \zeta\left ( \frac{\log 3}{\log2} \right ), 
\end{displaymath}
 and the
theorem follows.
\end{proof}
 
\begin{rem}\label{GI} In [GI2], Guido and Isola have associated a spectral triple to a
general self-similar fractal in $\br$.  For such a spectral triple, they proved
the same type of result as the one stated in Theorem \ref{proportional}. It is possible
that their proof of [GI2], Lemma 4.10, can be adapted to our present case,
from which Theorem \ref{proportional} would then follow by using the uniqueness of a
normalized self-similar measure on the gasket having the same homogeneity as
the $(\log 3/\log 2)$-Hausdorff measure. However, since the functional $\psi$
of Theorem \ref{psi} offers a nice description of the Hausdorff measure on the
Sierpinski gasket, we have preferred to give a direct argument which shows
that the positive functionals $\tau$ and $\psi$ on C$(\cs\cg)$  are
proportional.
\end{rem}

\begin{rem}\label{KL}
In the more delicate context of standard analysis on fractals,
Kigami and Lapidus identify in \cite{KiL2} the volume measure constructed in \cite{La3} via a
Dixmier trace functional. In particular, they show that this measure is
self-similar but is not always proportional to the natural Hausdorff measure
on the self-similar fractal, even when the measure maximizes the spectral
exponent (i.e., is an analogue of `Riemannian volume', in the sense of
\cite{La3,La4}). In the latter case, however, and for the special case of the
standard Sierpinski gasket, it does coincide with the natural Hausdorff
measure, pointing to some possible connections or analogies between the two
points of view.
\end{rem}

 \indent We will now discuss the concept of geodesic distance on $\cs\cg$ and
 then show that it can actually be measured by the metric induced by
 the spectral triple $ST(\cs\cg)$. The geodesic distance between two
 points $p$ and $q$  on the gasket is denoted $d_{geo}(p,q)$, and it 
 is defined as the minimal length of a rectifiable continuous curve on the gasket
 connecting $p$ and $q$. This metric on the gasket is studied in
 Barlow's lecture notes \cite{Ba} and from page 13 of
 these notes we quote the proposition below, which shows that the
 geodesic metric on the gasket is Lipschitz equivalent to the
 restriction of the Euclidean metric. In particular, the
 geodesic metric induces the standard topology on the gasket. 

\begin{Prop} \label{dgeqd}
Let $p, q$ be arbitrary points in $\cs\cg$ and let $\|p-q\|$ denote their
Euclidean distance. Then 
$$\|p-q\|\, \leq \, d_{geo}(p,q)\, \leq \, 8\|p - q\|.$$ 
\end{Prop}

\indent We have not found an easy way to give an analytic expression of the
geodesic distance between two arbitrary points in the gasket, 
but the distance from a corner point, a vertex, in $\cs\cg_0$ to an arbitrary
point in $\cs\cg$ can be 
precisely expressed in terms of barycentric coordinates, as we now show.

\begin{Lemma} \label{bacoor}
Let $v_1, v_2, v_3$ be the vertices  of $\cs\cg_0$. For any point $p$
in $\cs\cg$, let $(x, y, z) $ denote the barycentric coordinates
of $p$ with respect to $v_1, v_2, v_3$. Then the geodesic distance in
$\cs\cg$ from $v_1 $ to $p$ is $y+z$.
\end{Lemma}

\begin{proof}
Since the Euclidean metric is equivalent to the geodesic metric,
  and the function which assigns barycentric coordinates to
a point is continuous with respect to these metrics, we see
that it is enough to prove the statement for a point $p$ which is a
vertex in one of the triangles from one of the sets
$\cs\cg_n$ where $n$ is in $ \bn.$  
Let then $n$ in $\bn$ be fixed and let $p$ denote a vertex in $\cs\cg_n,$
and let us consider paths in $\cs\cg_n$ which connect $v_1$ to $p;$ we
will then show that a geodesic curve can be obtained inside $\cs\cg_n.$
Any path from $v_1$ to $p$ in $\cs\cg_n$  
will be a sum of steps, each of which will be
a positive multiple of one of the following 6  vectors, which all have
length $2 \pi / 3.$ 

$$
v_2-v_1, v_1- v_2, v_3 - v_1, v_1-v_3, 
v_3-v_2, v_2 - v_3.$$

When $p = (x,y,z) $ in barycentric coordinates, we have $p = v_1 + y(v_2-v_1)
+ z(v_3-v_1)$; so it must follow that the geodesic distance between
$v_1$ and $p$  is at least
$y+z$. On the other hand, it is possible, via a little drawing, to see
that $\cs\cg_n$ contains a path of length $y+z$  between $v_1$ and $p.
$ The geodesic distance between $v_1$ and $p$ is then $y+z,$ and this
result extends to a general point $p$ in $\cs\cg,$ by continuity.
\end{proof} 

\indent Let us then consider the geodesic 
 distance between two arbitrary, but different,  points $p$
and $q$. In order to describe a way to compute this distance, we will look
at the other picture of $\cs\cg,$ as the limit of a decreasing sequence,
say $\cf_n$, of compact subsets of the largest solid triangle, as
depicted in Figure \ref{SG2}.
 For any nonnegative integer  $n$, 
$\cf_n$ is the union of $3^n$ equilateral solid triangles, say $\{T_{n,k}\,|
1\leq k \leq 3^n\,\}$, with side length $2^{-n} 2 \pi / 3$. 
In order to determine
the geodesic distance between $p$ and $q$, we determine the largest
number $n_0 \in \bn_0$ for which  there exists a $k_0 $ in $\{1, 2,
\dots, 3^{n_0} \}$
such that both $p$ and $q$ belong to $T_{n_0,k_0}$. From the solid triangle
$T_{n_0,k_0}$ remain 3 solid triangles, say $T_{n_0+1,k_1},
\, T_{n_0+1,k_2}, \, T_{n_0+1,k_3}$, in the next step of
the iterative construction.
By assumption, $p$
and $q$ must lie in two different of these smaller triangles,
 say $T_{n_0+1,k_1}$ and
$T_{n_0+1,k_2}$. Hence, any path from $p$ to $ q$   must run from $p$
to a vertex in $T_{n_0+1,k_1}$. We can then use Lemma \ref{bacoor} to
measure the length of this part of the path. The vertex we have
arrived at may, or may not, be in $T_{n_0+1,k_2}$ too. In the first case,
we can measure the distance from the corner to $q$, again using
Lemma \ref{bacoor}. In the second case, the corner of  $T_{n_0+1,k_1}$
at which we have arrived must be directly connected to a corner of
$T_{n_0+1,k_2}$ via a side of the third triangle $T_{n_0+1,k_3}$. This
side has length $2^{-(n_0+1)} 2 \pi /3, $ and this edge will bring us
to a corner in $T_{n_0+1,k_2},$ from where we can measure the distance
to $q$ on the basis of Lemma \ref{bacoor}.
 This observation has several consequences, some of which we will formulate
in the next results.

\begin{Lemma} \label{dgeodiff}
Let $q$ be a point in $\cs\cg$ and let $g$ be the continuous function
on $\cs\cg$ defined by $g(p) = d_{geo}(p,q)$. Then, for any continuous
curve $r_{n,i} $ which parameterizes a triangle  in the gasket, we have 
$\|\,[D_{r_{n,i}}, \pi_{r_{n,i}}(g)]\,\| \leq 1.$ 
\end{Lemma}
\begin{proof}
Let us return to the introductory example in Lemma \ref{bacoor} and suppose
for simplicity that $q$ is the vertex $v_1$.
 Given an edge, say $e$, 
 in a  triangle $\Delta_{n,i} $, which is parameterized by one of the 
 functions $r_{n,i},$       
then there are  only  3 possibilities for the slope of the  edge,
since  it must be parallel to one of the edges of the big
triangle. Following Lemma \ref{bacoor}, we therefore deduce that in this
case where $q = v_1,$  we must have 
\begin{align*}
 \exists \beta \in \br, \exists \alpha \in \{-1,
0, 1\}:& \\  r_{n,i}(t) \in e \, \Rightarrow \,&  g(r_{n,i}(t))\, = \,
d_{geo}(r_{n,i}(t),v_1)  \, = \, \alpha t + \beta. 
\end{align*}
According to Lemma  \ref{domD}, such a function is in the
domain of all the operators  $D_{n,i},$ and we see that the
derivative is numerically bounded by $1$. We will now establish a similar
result for a general point $q.$ This situation is discussed in the
text just in front of this lemma and it follows from there that the function
$g(r_{n,i}(t))$ has to be modified by an additive constant, which
measures the geodesic distance from $q$ to one of the endpoints of the
edge $e.$ When $r_{n,i}(t)$ passes along the edge $e,$ the
geodesic distance from $q$ may reach an extremal value in the interior
of the edge, and the slope
may change, but still be in the set $\{-1, 0, 1\}.$ 
 Hence,  according to Lemma  \ref{domD}, 
we deduce that the derivative of $g(r_{n,i}(t))$ exists
almost everywhere and is numerically bounded by 1, and by the same
lemma we get $\|\,[D_{r_{n,i}}, \pi_{r_{n,i}}(g)]\,\| \leq 1.$
\end{proof}
 \begin{Thm} \label{Sgthm}
The metric on the Sierpinski gasket, $d_{\cs\cg},$ induced
by the spectral triple $(\mathrm{C}(\cs\cg), H_{\cs\cg}, D_{\cs\cg}),$
coincides with the geodesic distance on $\cs\cg.$
\end{Thm}

\begin{proof}
As before, let $\cn = \{ f \in \mathrm{C}(\cs\cg)\,|\, \|\,[D_{\cs\cg},
\pi_{\cs\cg} (f)]\,\| \leq 1 \,\}$. Then Lemma \ref{dgeodiff}
 shows that  all
the functions of the form $g(p) = d_{geo}(p,q) $ belong to $\cn$. We
therefore have 
$$ \forall p, q \in \cs\cg: \quad d_{geo}(p,q) \leq d_{\cs\cg}(p,q).$$
The other inequality is an application of the fundamental theorem of
calculus. Let points $p, q $ in $\cs\cg$ be given and let $f$ be in   $
\cn.$ Then $f$ is continuous and for
 any $\eps > 0$ there exists a natural number $n$ and vertices
$v, w $ in $\cs\cg_n$ such that 
\begin{equation} \label{epsi}
d_{geo}(p,v) + d_{geo}(q,w) + |f(p) - f(v)| + |f(q) - f(w)| < \eps.
\end{equation}  
 Let us then look at
$|f(v) - f(w)|$  and  show that this quantity is at most $d_{geo}(v,w).$
We will therefore consider a path along some of the edges
of some of the triangles $\Delta_{n,i}$
 connecting $v$ and $w$ inside $\cs\cg_n.$ 
Since $ f$ is in $\cn$, it follows from
 Lemma \ref{domD}  that the derivative of the induced function 
$f(r_{n,i}(t))$
corresponding to    an edge in $\Delta_{n,i} $
is a measurable function which is numerically  bounded
by 1 almost everywhere. We can then express $f(v) - f(w) $ as a sum of line
integrals. Each of these integrals is numerically bounded from above by the length
of the path over which the integration is performed; so $|f(v) - f(w) |$
is dominated by the  length  of any path in $\cs\cg_n$ connecting $s$ and $t$.
This implies that

\begin{displaymath} 
   |f(v) - f(w)|  \, \leq \, d_{geo}(v,w) \quad \, \, 
\end{displaymath}
\noindent
 and, by equation (\ref{epsi}),
\begin{displaymath}
    |f(p) - f(q)|  \,  \leq \, d_{geo}(p,q) + \eps.
\end{displaymath}

\noindent
From here we see that 
$$d_{\cs\cg}(p,q ) \, = \, \underset{f \in \cn}{\sup}|f(p) - f(q)| \,
\leq \, d_{geo}(p,q),$$
and the theorem follows.
\end{proof}

\section{Concluding remarks}

We close this paper by indicating several possible extensions 
or directions of future research related to this work, some of which may
be investigated in later articles:

\smallskip

\begin{itemize}

\item[(i)] exploring the possibility of associating to other 
fractals built on curves spectral triples which 
are based on a direct sum of $r$-tri-\newline ples  
and describing the topological and geometric properties of those fractals; 
\smallskip
\item[(ii)] investigating the geometric and topological properties of the Sierpisnki gasket described  
by the spectral triples based on some other possible direct sums of spectral triples for circles;  
\smallskip
\item[(iii)] extending connections and/or relations to other known constructions of
differential operators on the Sierpinski gasket and other
self-similar fractals, especially in some of the approaches to 
{\em analysis on fractals} expounded, for example, in \cite{Ba} and
\cite{Ki}, either probabilistically or analytically. In particular, 
investigating some of the open problems or conjectures proposed 
within that framework in \cite{La3,La4}, and aimed at 
merging aspects of fractal, spectral and noncommutative geometry; 
\smallskip
\item[(iv)] studying the differential operators
(including `Laplacians') connected to the Dirac-type
operators constructed in this paper, 
as well as of the solutions 
of partial differential equations naturally associated
to them;
\smallskip
\item[(v)] looking at the Sierpinski gasket
via the `harmonic coordinates' attached to the
Laplacian $\Delta$ associated to our Dirac operator $D$
(namely, $-\Delta = D^2$), or to a suitable modification thereof. 
 (See e.g. \cite{Te1, Te2}
and the relevant references therein for the analogous
situation involving the usual Laplacian on the gasket.);
\smallskip
 \item[(vi)] looking for further applications of the spectral triples as a tool for
  computing invariants of algebraic topological type for the
  Sierpinski gasket and other fractals.
\end{itemize}
\smallskip

\indent We now explain more precisely what we mean by item (ii): 
\smallskip

\indent We first recall the construction where the Sierpinski gasket 
is obtained as the limit of an increasing sequence of sets in $\br^2$ (see Figure 1). 
We take as a starting point an equilateral triangle with side length $2 \pi / 3$ and we call it $\cs\cg_0$. 
The next figure, $\cs\cg_1$, is obtained by adding another triangle of half size,
and turned upside down relative to  $\cs\cg_0$, and so on. 
This means that $\cs\cg_1$ consists of  $\cs\cg_0$,  and 1 equilateral 
triangle of size $2^{-1} $ of $\cs\cg_0$.  The following figure $\cs\cg_2$ 
then consists of $\cs\cg_0$, 1 equilateral triangle of size $2^{-1} $ of 
$\cs\cg_0$, and 3 of size  $2^{-2}$ of $\cs\cg_0$. And finally, $\cs\cg_n$ 
consists of $\cs\cg_0$, 1 equilateral triangle of size $2^{-1} $ of $\cs\cg_0$, 
3 of size  $2^{-2}$ of $\cs\cg_0$, and so on up to $3^{n-1}$ triangles of size $2^{-n}$ that of the starting one.   
We will use  a circle as the basis for the construction of
a spectral triple for each of the triangles contained in  $\cs\cg_n$ as $n\rightarrow \infty$, and
then ultimately  obtain a spectral triple for the Sierpinski gasket in the 
same way as it was done in Section 8 of the present article.  Another possibility is to use 
a circle as the basis for the construction of
a spectral triple for each of the triangles considered in Section 8, as well as each of the triangles 
considered just above, 
and then ultimately  obtain a spectral triple for the Sierpinski gasket 
in the same way as it was done in Section 8.
\smallskip

\indent Next, in the special case of the Sierpinski gasket, we
elaborate on item (vi) above regarding the computation of topological invariants:

\indent The spectral triple $ST(\cs\cg)$ has been constructed in such a way 
that the theory of \cite{Co2}, Sections IV.1 and IV.2, will yield nontrivial 
topological information. The details of this will appear in a
later article, but to indicate what sort of results we have in mind,
we will just describe the bounded Fredholm module (K-cycle) which one can
obtain from $ST(\cs\cg) = (C(\cs\cg), H_{\cs\cg}(\text{together with } \pi_{\cs\cg}), D_{\cs\cg}).$
The unitary part, $F_{\cs\cg} $, of the polar decomposition of the Dirac
operator $D_{\cs\cg}$ gives a bounded odd Fredholm module
$(\pi_{\cs\cg}, H_{\cs\cg}, F_{\cs\cg}).$ The pairing with
$K_1(\cs\cg)$ can then be obtained in the following way. Let
$F_{\cs\cg} = P_+ -(I-P_+)$, where $P_+$ is the orthogonal projection
of $H_{\cs\cg}$  onto the eigenspace 
corresponding to the positive eigenvalues of $D_{\cs\cg},$ and let $f$ be a continuous function on
$\cs\cg $ which has the property that $\pi_{\cs\cg}(f)$ is invertible,
i.e. $\min\{ |f(x)| \, |\, x \in \cs\cg \} \, > \, 0. $ Then the
equicontinuity of $f$ will imply that there exists a natural number,
say $n_0$, such that for any $n \geq n_0$ and any triangle
$\Delta_{n,i}$ of size $2^{-n}$ that of the original one, the function $f$
has winding number 0 along $\Delta_{n,i}$. This implies that the
operator $P_+ \pi_{\cs\cg}(f)P_+$ has a finite index which
actually is the opposite number to the sum of all the winding numbers
over all the triangles $\Delta_{j,h} $ for $0 \leq j < n_0$ and $h \in
\{1, \dots 3^j\}$. This is a rather obvious extension of well-known
results for the circle, but we think that also the newer
invariants related to the cyclic cohomology of the gasket may be
expressed via the spectral triple $ST(\cs\cg).$


\begin{thebibliography}{99}


\bibitem[1]{Ba}
M. T. Barlow, Diffusions on fractals, in: Lectures on
Probability Theory and Statistics (Saint-Flour, 1995),  
 Springer Lecture Notes 1690 (1998)  1--121.  


\bibitem[2]{CI1}
E. Christensen, C. Ivan,   Sums of two dimensional spectral
  triples, Math. Scand. 100 (2007) 35--60.
  
\bibitem[3]{CI2}
E. Christensen, C. Ivan, Spectral triples for AF
  C*-algebras and metrics on the Cantor set, J. Operator Theory
 56 (2006) 17--46.
  
\bibitem[4]{CiGS}
F. Cipriani, D. Guido, S. Scarlatti, A remark on trace properties 
of \newline 
K-cycles, J. Operator Theory 35 (1996) 175--189.


\bibitem[5]{Co1}
A. Connes, Compact metric spaces, Fredholm modules, and
  hyperfiniteness,  Ergodic. Theor. and Dynam. Systems 9 (1989) 207--220.
    
\bibitem[6]{Co2} A. Connes, Noncommutative Geometry, Academic
    Press, San Diego (1994).
    
\bibitem[7]{Co3}
A. Connes, Unpublished notes on a Dirac operator associated to
  the Cantor subset of the unit interval. (Electronic message to Michel Lapidus, May 2002.)
  
\bibitem[8]{CoM}
A. Connes, M. Marcolli, A walk in the noncommutative garden, arXiv.math.QA/0601054 (2006). 

\bibitem[9]{CoS}
A. Connes, D. Sullivan, Quantized calculus on $S^1$ and quasi-Fuchsian
groups,  unpublished (1994)..


\bibitem[10]{CMRV}
G. Cornelissen, M. Marcolli, K. Reihani, A. Vdovina,  Noncommutative geometry on trees and buildings,
arXiv.math.QA/0604114  (2006).


\bibitem[11]{Da}
E. B. Davies, Analysis on graphs and noncommutative geometry,  J. Funct. Anal. 111 (1993) 398-­430. 


\bibitem[12]{Edg1}
G. A. Edgar, Integral, Probability, and Fractal Measures, Springer-Verlag, New York (1998).


\bibitem[13]{Edg2}
G. A. Edgar, Measure, Topology, and Fractal Geometry, Undergraduate Texts in Mathematics,
Springer-Verlag, New York (1990). 


\bibitem[14]{Edw}
R. E. Edwards, Fourier Series. A modern introduction, Vol. 1, Second Edition, Graduate Texts in Mathematics 64, 
Springer-Verlag, New York (1979). 


\bibitem[15]{Fa}
K. J. Falconer,  Fractal Geometry: Mathematical foundations and
applications, Wiley, Chichester (1990). 


\bibitem[16]{GI1} D. Guido, T. Isola, Fractals in noncommutative geometry,
in: Proc. Conf. ``Mathematical Physics in Mathematics and Physics"
Siena 2000, ed. R. Longo, Fields Institute Communications,
Amer. Math. Soc. 30 (2001) 171--186.


\bibitem[17]{GI2} D. Guido, T. Isola, Dimensions and singular
      traces for spectral triples, with applications to fractals,
    J.  Funct.  Anal.  203 (2003) 362--400.


\bibitem[18]{GI3} D. Guido, T. Isola,  Dimensions and spectral
    triples for fractals in $\br^n$, arXiv.math.OA/0404295 v2 (2005).
    

\bibitem[19]{KaR}
R. V. Kadison, J. R. Ringrose,  Fundamentals of the Theory of
  Operator Algebras, Vols. I and II, Academic Press, New York (1983).


\bibitem[20]{Ki} J. Kigami, Analysis on Fractals, 
Cambridge Univ. Press, Cambridge  (2001). 

  
\bibitem[21]{KiL1} J. Kigami, M. L. Lapidus, Weyl's problem for the spectral
distribution of Laplacians on p.c.f. self-similar sets, Commun. Math. Phys. 158 (1993)
  93--125.
  
  
 \bibitem[22]{KiL2} J. Kigami, M. L. Lapidus, Self-similarity of volume measures
for Laplacians on p.c.f. self-similar fractals, Commun. Math. Phys. 217 (2001)
  165--180.  
  
  
 \bibitem[23]{Ku1} P. Kuchment, Quantum graphs I.
Some basic structures, Waves in Random Media 14 (2004)
  107--128.  
  
  
 \bibitem[24]{Ku2} P. Kuchment,  Quantum graphs II.
Some spectral properties of quantum and combinatorial graphs, 
J. Phys. A  38 (2005) 4887--4900.  


\bibitem[25]{La1} M. L. Lapidus,  Fractal drum, inverse spectral
    problems for elliptic operators and a partial resolution of the
    Weyl--Berry conjecture, Trans. Amer. Math. Soc. 325 (1991)
  465--529.  
  

\bibitem[26]{La2} M. L. Lapidus, Vibrations of fractal drums,
    the Riemann hypothesis, waves in fractal media, and the Weyl--Berry
    conjecture, in: Ordinary and Partial Differential Equations, eds. B.
  D. Sleeman and R. J. Jarvis, Pitmann Research Notes in Math. 325
  (1993) 126--209.
  
  
\bibitem[27]{La3}  M. L. Lapidus,  Analysis on fractals,
      Laplacians on self-similar sets, noncommutative geometry and
      spectral dimensions, Topological Methods in Nonlinear Analysis
   4 (1994) 137--195.  


\bibitem[28]{La4} M. L. Lapidus,  Towards a noncommutative
    fractal geometry? Laplacians and volume measures on fractals, in:
  Harmonic Analysis and Nonlinear Differential Equations,
  Contemp. Math. Amer. Math. Soc.  208 (1997) 211--252.



\bibitem[29]{LM} M. L. Lapidus, H. Maier,  The Riemann hypothesis
    and inverse spectral problems for fractal strings, J. London
  Math. Soc. (2)  52 (1995) 15--34.
  
\bibitem[30]{LP} M. L. Lapidus, C. Pomerance,  The Riemann
    zeta-function and the one-dimensional Weyl--Berry conjecture for
    fractal drums, Proc. London Math. Soc. (3) 66 (1993) 41--69.  
    
    
\bibitem[31]{L-vF1} M. L. Lapidus, M. van Frankenhuysen,  Fractal
    Geometry and Number Theory: Complex dimensions of fractals strings 
and zeros of zeta functions, Birkh\" auser, Boston (2000). 


\bibitem[32]{L-vF2} M. L. Lapidus, M. van Frankenhuijsen,  Fractal
    Geometry, Complex Dimensions and Zeta Functions:
Geometry and spectra of fractals strings, Springer Monographs in Mathematics,
Springer-Verlag, New York (2006).



\bibitem[33]{Re} M. Requardt, Dirac operators and the calculation of 
the Connes metric on arbitrary (infinite) graphs, J. Phys. A 35 (2002) 759--779. 


\bibitem[34]{Ri1}
M. A. Rieffel, Comments concerning non-commutative metrics, 
Talk given at an AMS Special Session, Texas A$\&$M (1993).


\bibitem[35]{Ri2}
M. A. Rieffel,  Metrics on states from actions of compact groups,
Doc. Math. 3 (1998) 215--229.

\bibitem[36]{Ri3}
M. A. Rieffel,  Metrics on state spaces,
Doc. Math. 4 (1999) 559--600.

\bibitem[37]{Ri4}
M. A. Rieffel,  Compact quantum metric spaces, in: 
Operator Algebras, Quantization, and Non-Commutative Geometry,
Contemp. Math. Amer. Math. Soc. 365 (2004) 315--330.


\bibitem[38]{SWW} 
E. Schrohe, M. Walze, J. M. Warzecha,
 Construction de triplets spectraux \`a partir de modules de
 Fredholm,  C. R. Acad.  Sci.  Paris S\'{e}r.  I Math.  326 (1998) 1195--1199.
    

\bibitem[39]{St} 
R. S. Strichartz,  Analysis on fractals, Notices
Amer. Math. Soc. (10) 46 (1999)  1199--1208. 



\bibitem[40]{Te1} 
A. Teplyaev, Energy and Laplacian on the Sierpinski gasket, in:
Fractal Geometry and Applications: A jubilee of Beno\^it Mandelbrot,
eds. M. L. Lapidus and M. van Frankenhuijsen, 
Proc. Sympos. Pure Math. Amer. Math. Soc. 72 (2004)  131--154. 


\bibitem[41]{Te2}
A. Teplyaev, Harmonic coordinates on fractals with finitely ramified cell structure,
arXiv:math.PR/0506261 (2006).    

    
\end{thebibliography}
\end{document}